\documentclass[mathpazo]{cicp}
\usepackage{xcolor}
\usepackage{fancyhdr}
\usepackage{graphicx}
\usepackage{subfigure}
\usepackage{amssymb}
\usepackage{amsmath}
\usepackage{amsfonts}
\usepackage{theorem}
\usepackage{mathrsfs}
\usepackage{mathtools}
\usepackage{bm}
\usepackage{color}
\usepackage{amsthm}
\usepackage{setspace}
\usepackage{exscale}
\usepackage{relsize}
\usepackage{epstopdf}
\usepackage{float}
\usepackage{cite}
\usepackage{color}
\usepackage{nicefrac}

\usepackage{booktabs,multirow} 
\usepackage{array,amsmath} 
\usepackage{paralist} 
\usepackage{verbatim} 
\usepackage{subfigure}

\newcommand*\xbar[1]{%
  \hbox{%
    \vbox{%
      \hrule height 0.5pt 
      \kern0.4ex
      \hbox{%
        \kern-0.05em
        \ensuremath{#1}%
        \kern-0.00em
      }%
    }%
  }%
}
\newcommand\eref[1]{(\ref{#1})}

\newcommand{\mF}{\bm{F}}

\newcommand{\mG}{\bm{G}}

\newcommand{\mU}{\bm{U}}

\newcommand{\dt}{\Delta t}
\newcommand{\dx}{\Delta x}
\newcommand{\dy}{\Delta y}
\newcommand{\eps}{\varepsilon}

\newcommand{\hf}{{\frac{1}{2}}}

\newcommand{\jph}{{j+\frac{1}{2}}}
\newcommand{\jmh}{{j-\frac{1}{2}}}
\newcommand{\kph}{{k+\frac{1}{2}}}
\newcommand{\kmh}{{k-\frac{1}{2}}}

\newcommand{\ajphp}{{a_{j+\frac{1}{2}}^+}}
\newcommand{\ajphm}{{a_{j+\frac{1}{2}}^-}}

\newcommand{\bmF}{\bm{\mathcal{F}}}
\newcommand{\bmG}{\bm{\mathcal{G}}}
\begin{document}
\title[ovel Adaptive Schemes]{Novel Adaptive Schemes for Hyperbolic Conservation Laws}
\author[Chu S et.~al.]{Shaoshuai Chu\affil{1}, Pingyao Feng\affil{2}, Vadim A. Kolotilov\affil{3}, Alexander Kurganov\affil{4}\comma\corrauth, and Vladimir V. Ostapenko\affil{3}}
\address{\affilnum{1} Department of Mathematics, RWTH Aachen University, 52056 Aachen, Germany.\\
\affilnum{2} Department of Mathematics, North Carolina State University, 27695 Raleigh, USA.\\
\affilnum{3} Lavrentyev Institute of Hydrodynamics, Siberian Branch, Russian Academy of Sciences, Novosibirsk, 630090 Russia. \\
\affilnum{4} Department of Mathematics and ShenZhen International Center for Mathematics, Southern University of Science and Technology, Shenzhen, 518055, P.R. China.\\
}
\email{{\tt chu@igpm.rwth-aachen.de} (S.~Chu),  {\tt pfeng3@ncsu.edu} (P.~Feng),  {\tt kolotilov1992@gmail.com} (V. A.~ Kolotilov), ~ {\tt alexander@sustech.edu.cn} (A.~Kurganov),  {\tt ostigil@mail.ru} (V. V. Ostapenko).}


\begin{abstract}
We introduce new adaptive schemes for the one- and two-dimensional hyperbolic systems of conservation laws. Our schemes are based on an
adaption strategy recently introduced in [{\sc S. Chu, A. Kurganov, and I. Menshov}, Appl. Numer. Math., 209 (2025)]. As there, we use a
smoothness indicator (SI) to automatically detect ``rough'' parts of the solution and employ in those areas the second-order finite-volume
low-dissipation central-upwind scheme with an overcompressive limiter, which helps to sharply resolve nonlinear shock waves and linearly
degenerate contact discontinuities. In smooth parts, we replace the limited second-order scheme with a quasi-linear fifth-order (in space
and third-order in time) finite-difference scheme, recently proposed in [{\sc V. A. Kolotilov, V. V. Ostapenko, and N. A. Khandeeva},
Comput. Math. Math. Phys., 65 (2025)]. However, direct application of this scheme may generate spurious oscillations near ``rough'' parts,
while excessive use of the overcompressive limiter may cause staircase-like nonphysical structures in smooth areas. To address these issues,
we employ the same SI to distinguish contact discontinuities, treated with the overcompressive limiter, from other ``rough'' regions, where
we switch to the dissipative Minmod2 limiter. Advantages of the resulting adaptive schemes are clearly demonstrated on a number of
challenging  numerical examples. 
\end{abstract}


\maketitle

\section{Introduction}
\label{sec1}
We consider hyperbolic systems of conservation laws, which in the one-dimensional (1-D) and two-dimensional (2-D) cases, read as
\begin{equation}
\mU_t+\mF(\mU)_x=\bm0,
\label{1.1}
\end{equation}
and
\begin{equation}
\mU_t+\mF(\mU)_x+\mG(\mU)_y=\bm0,
\label{1.2}
\end{equation}
respectively. Here, $x$ and $y$ are spatial variables, $t$ is the time, $\mU\in\mathbb R^d$ is a vector of unknown functions, and
$\mF:\mathbb R^d\to\mathbb R^d$ and $\mG:\mathbb R^d\to\mathbb R^d$ are fluxes.

It is well-known that solutions of \eref{1.1} and \eref{1.2} can produce extremely complex wave structures including shocks, rarefactions,
and contact discontinuities even when the initial data are smooth. This makes it quite challenging to develop accurate and reliable
numerical methods for \eref{1.1} and \eref{1.2}. Such methods are typically based on equivalent integral forms of \eref{1.1} and \eref{1.2}
for the time evolution of the computed solution and on nonlinear limiters, designed to prevent spurious oscillations. For a variety of
existing numerical methods for hyperbolic systems of conservation laws, we refer the reader to the monographs and review papers
\cite{GR2,Tor,Leveque02,Hesthaven18,BAF,Shu20} and references therein.

In this paper, we focus on scheme adaption methods, which are based on automatic detection of smooth and nonsmooth (``rough'') parts of the
computed solutions and on applying different numerical methods in the corresponding parts of the computational domain. The main advantage of
such an approach is related to the fact that special nonlinear limiting techniques can be used to accurately capture shock waves and contact
discontinuities only, while smooth parts of the solution may be accurately and efficiently computed by much simpler, possibly higher-order,
(quasi-linear) unlimited methods.

A success of adaptive schemes hinges on an accuracy and robustness of the smoothness indicator (SI) used to detect ``rough'' parts of the
computed solution. A wide variety of different SIs has been developed; see, e.g., 
\cite{DZLD14,Dewar15,ABD08,GLHWDW,GT06,GT02,GPP,PupSem,QS05,VR16,FS17,WSYK15,KKP02,RL87,WDGH20,CCK23_Adaptive}. Another critical component
of a successful adaptive scheme is a very sharp method used to capture shock and contact discontinuities inside the identified ``rough''
parts of the solution. In this paper, we follow \cite{CKM2025} and implement the SI proposed in \cite{RL87} along with the second-order
semi-discrete finite-volume low-dissipation central-upwind (LDCU) schemes from \cite{CKX_24} combined with the overcompressive and
dissipative limiters introduced in \cite{Lie03}.

The novelty of the proposed adaptive schemes is twofold. First, we treat the smooth parts of the computed solutions in a different way.
Instead of using the second-order LDCU scheme with the dissipative Minmod2 limiter (see \cite{Lie03} for the classification of different
limiters) as it was done in \cite{CKM2025}, we increase the resolution by switching to the recently proposed quasi-linear fifth-order (in
space and third-order in time) finite-difference scheme recently introduced in \cite{KOK25}. Following \cite{LNOT23}, we call a numerical
scheme {\em quasi-linear} if it becomes linear when applied to linear hyperbolic systems. As it was shown in \cite{KOK25}, the fifth-order
quasi-linear scheme is very accurate since it is fifth-order accurate in space in all smooth parts of the computed solution,
including those smooth subregions that intersect with the shock influence domains. However, the fifth-order quasi-linear scheme does not
contain any nonlinear stabilization mechanism and thus, as we show in \S\ref{sec4}, it can produce spurious oscillations near
the interfaces between the ``rough'' and ``smooth'' regions. To suppress these oscillations, one can ``enlarge'' the identified ``rough''
areas, but in this case, the overcompressive limiter would affect the smooth parts of the solution causing development of artificial kinks
and staircase-like structures or even nonphysical jump discontinuities there. 

In order to address these issues, we introduce an alternative approach, which is the second novel component of the proposed adaptive
schemes. To this end, we use the same SI, but applied to different fields to automatically distinguish between the contact waves and other
``rough'' parts of the solution. For example, in the case of Euler equations of gas dynamics, we take advantage of the fact that while the
density jumps at the contact surfaces, the pressure remains continuous there. After identifying the locations of contact surfaces, we apply
the second-order LDCU scheme with an overcompressive limiter in their vicinities only, while utilizing the dissipative Minmod2 limiter in
the rest of the ``rough'' areas. In addition, the overcompressive and Minmod2 limiters are applied to the local characteristic variables
using the local characteristic decomposition (LCD), which is often used in the context of high-order schemes \cite{Qiu02} to suppress
spurious oscillations, but can also be implemented to enhance the resolution of second-order schemes; see, e.g., \cite{CKX24,CHK25}. This
way, we construct efficient and robust adaptive schemes, which achieve a very high resolution of both jump discontinuities and smooth parts
of the computed solutions.

The rest of the paper is organized as follows. In \S\ref{sec2}, we first review the 1-D LDCU scheme from \cite{CKX_24} and the
quasi-linear fifth-order finite-difference scheme from \cite{KOK25}, and then introduce a new 1-D adaptive scheme,
which is based on the dynamical splitting of the computational domains into three parts: vicinities of the contact discontinuities, other
``rough'' parts of the computed solutions, and smooth parts of the computed solutions. In \S\ref{sec3}, we extend the new adaptive scheme
to the 2-D case. In \S\ref{sec4}, we test the developed schemes on a number of 1-D and 2-D numerical examples for the Euler equations of gas
dynamics. The obtained numerical results demonstrate that the adaptive schemes contain a substantially smaller amount of numerical
dissipation and achieve higher resolution compared with the adaptive schemes from \cite{CKM2025}. Finally, in \S\ref{sec6}, we provide
concluding remarks.

\section{Novel 1-D Adaptive Scheme}\label{sec2}
In this section, we consider the 1-D hyperbolic system \eref{1.1} and describe a scheme adaption algorithm similar to the one recently
introduced in \cite{CKM2025}. The novel feature here is in the treatment of smooth parts of the computed areas detected by the SI; see
\S\ref{sec22} and \S\ref{sec23}.

\medskip 
\subsection{1-D Low-Dissipation Central-Upwind (LDCU) Schemes}\label{sec21}
In this section, we briefly review the second-order semi-discrete LDCU schemes recently introduced in \cite{CKX_24}.

Assuming that the computational domain is covered with the uniform cells $I_j:=\big[x_\jmh,x_\jph\big]$ with $x_\jph-x_\jmh\equiv\dx$
centered at $x_j=\big(x_\jmh+x_\jph\big)/2$, we suppose that the computed cell averages,
\begin{equation*}
\xbar\mU_j(t):\approx\frac{1}{\dx}\int\limits_{I_j}\mU(x,t)\,{\rm d}x,
\end{equation*}
are available at a certain time level $t\ge0$. Note that, $\,\xbar\mU_j(t)$ like many other indexed quantities, which will be introduced
below, are time-dependent, but from here on, we suppress their time-dependence for the sake of brevity.

The computed cell averages are evolved in time by numerically solving the following system of ODEs:
\begin{equation}
\frac{{\rm d}\xbar\mU_j}{{\rm d}t}=-\frac{\bm{{\cal F}}^{\rm LDCU}_\jph[\bm U]-\bm{{\cal F}}^{\rm LDCU}_\jmh[\bm U]}{\dx},
\label{2.1}
\end{equation}
where $\bm{{\cal F}}^{\rm LDCU}_\jph[\bm U]$ are the LDCU numerical fluxes given by (see \cite{CKX_24})
\begin{equation*}
\bm{{\cal F}}^{\rm LDCU}_\jph[\bm U]=\frac{\ajphp\mF\big(\mU^-_\jph\big)-\ajphm \mF\big(\mU^+_\jph\big)}{\ajphp-\ajphm}+
\frac{\ajphp\ajphm}{\ajphp-\ajphm}\Big(\mU^+_\jph-\mU^-_\jph\Big)+\bm q_\jph.
\end{equation*}
Here:

\smallskip
\noindent
$\bullet$ $\mU^\pm_\jph$ are the left- and right-sided point values of $\mU$ at the cell interface $x=x_\jph$ reconstructed out of the given
set of cell averages $\{\,\xbar\mU_j\}$ using a proper nonlinear limiter; see Appendix \ref{sec211};

\smallskip
\noindent
$\bullet$ $a^\pm_\jph$ are the one-sided local speeds of propagation, which can be estimated using the largest and smallest eigenvalues of
the Jacobian $A(\mU):=\frac{\partial\mF}{\partial\mU}(\mU)$;

\smallskip
\noindent
$\bullet$ $\bm q_\jph$ is a ``built-in'' anti-diffusion term, which needs to be derived for a particular system at hand.

\smallskip
In this paper, we consider the Euler equations of gas dynamics, which, in the 1-D case, read as \eref{1.1} with
\begin{equation}
\mU=\big(\rho,\rho u,E\big)^\top\quad{\rm and}\quad\mF=\big(\rho u,\rho u^2+p,u(E+p)\big)^\top.
\label{2.2}
\end{equation}
Here, $\rho$, $u$, $p$, and $E$ are the density, velocity, pressure, and total energy, respectively, and the system is closed through the
following equations of state (EOS) for ideal gases:
\begin{equation}
p=(\gamma-1)\Big[E-\hf\rho u^2\Big],
\label{2.3}
\end{equation}
where the parameter $\gamma$ represents the specific heat ratio. For the Euler system \eref{1.1},\\ \eref{2.2}--\eref{2.3}, the
one-sided local speeds of propagation can be estimated by
\begin{equation*}
a^+_\jph=\max\Big\{u^+_\jph+c^+_\jph,u^-_\jph+c^-_\jph,0\Big\},\quad a^-_\jph=\min\Big\{u^+_\jph-c^+_\jph,u^-_\jph-c^-_\jph,0\Big\},
\end{equation*}
where
$$
u^\pm_\jph=\frac{(\rho u)_\jph^\pm}{\rho^\pm_\jph},\quad p^\pm_\jph=(\gamma-1)\Big[E^\pm_\jph-\hf\rho^\pm_\jph\big(u^\pm_\jph\big)^2\Big],
\quad c^\pm_\jph=\sqrt{\frac{\gamma p^\pm_\jph}{\rho^\pm_\jph}},
$$
and the anti-diffusion term $\bm q_\jph$ is given by (see \cite{CKX_24})
\begin{equation*}
\bm q_\jph=\alpha^*_\jph q^\rho_\jph\left(1,u^*_\jph,\hf\big(u^*_\jph\big)^2\right)^\top.
\end{equation*}
Here,
\begin{equation*}
q^\rho_\jph={\rm minmod}\Big(\big(u^*_\jph-a^-_\jph\big)\big(\rho^*_\jph-\rho^-_\jph\big),
\big(a^+_\jph-u^*_\jph\big)\big(\rho^+_\jph-\rho^*_\jph\big)\Big),
\end{equation*}
$u^*_\jph=(\rho u)^*_\jph/\rho^*_\jph$ with $\rho^*_\jph$ and $(\rho u)^*_\jph$ being the first and second components of
\begin{equation*}
\mU^*_\jph=\frac{a^+_\jph\mU^+_\jph-a^-_\jph\mU^-_\jph-\Big[\bm F\big(\mU^+_\jph\big)-\bm F\big(\mU^-_\jph\big)\Big]}{a^+_\jph-a^-_\jph},
\end{equation*}
\begin{equation*}
\alpha^*_\jph=\left\{\begin{aligned}
&\frac{a^+_\jph}{a^+_\jph-u^*_\jph}&&\mbox{if}~u^*_\jph<0,\\&\frac{a^-_\jph}{a^-_\jph-u^*_\jph}&&\mbox{otherwise},\end{aligned}\right.
\end{equation*}
and the minmod function is defined by ${\rm minmod}(a,b):=\frac{{\rm sgn}(a)+{\rm sgn}(b)}{2}\min(|a|,|b|)$.

\noindent{\bf Remark 2.1\,\,}
The ODE system \eref{2.1} has to be numerically integrated using a sufficiently accurate and stable ODE solver. In the numerical examples
reported in this paper, we have used the three-stage third-order strong stability preserving (SSP) Runge-Kutta method, which reads as (see
\cite{Gottlieb01,Gottlieb11})
\begin{equation}
\begin{aligned}
&\mU^{\rm I}_j=\mU_j(t)-\frac{\dt}{\dx}\left(\bm{{\cal F}}_\jph-\bm{{\cal F}}_\jmh\right),\\
&\mU^{\rm II}_j=\frac{3}{4}\,\mU_j(t)+\frac{1}{4}\left[\mU^{\rm I}_j-
\frac{\dt}{\dx}\left(\bm{{\cal F}}^{\rm I}_\jph-\bm{{\cal F}}^{\rm I}_\jmh\right)\right],\\
&\mU_j(t+\dt)=\frac{1}{3}\,\mU_j(t)+\frac{2}{3}\left[\mU^{\rm II}_j-
\frac{\dt}{\dx}\left(\bm{{\cal F}}^{\rm II}_\jph-\bm{{\cal F}}^{\rm II}_\jmh\right)\right],
\end{aligned}
\label{2.8}
\end{equation}
where $\bm{{\cal F}}_\jph=\bm{{\cal F}}^{\rm LDCU}_\jph[\mU]$, $\bm{{\cal F}}^{\rm I}_\jph=\bm{{\cal F}}^{\rm LDCU}_\jph[\mU^{\rm I}]$, and $\bm{{\cal F}}^{\rm II}_\jph=\bm{{\cal F}}^{\rm LDCU}_\jph[\mU^{\rm II}]$. 
In \eref{2.8}, the time step $\dt$ is selected adaptively based on the following CFL-based stability restriction:
\begin{equation*}
\dt\le\frac{\dx}{2a},\quad a:=\max_j\Big\{\max\big(a^+_\jph,-a^-_\jph\big)\Big\}.
\end{equation*}

Notice that in \eref{2.8} we have replaced the cell averages of $\bm U$ with its point values $\bm U_j(t):\approx\bm U(x_j,t)$. This
does not affect the accuracy of the second-order scheme as the difference between the cell averages and point values is proportional to
${\cal O}((\dx)^2)$ for smooth solutions. On the other hand, switching to the finite-difference evolution is advantageous as the LDCU scheme
is going to be combined with the quasi-linear fifth-order finite-difference scheme from \cite{KOK25}, for which the point values of the
solution are to be evolved in time.

\subsection{1-D Quasi-Linear Fifth-Order Scheme}\label{sec22}
In this section, we describe the quasi-linear fifth-order finite-difference scheme recently introduced in \cite{KOK25}. As this scheme is
going to be implemented in our adaption algorithm, we write it in the fully discrete flux form \eref{2.8}.

Assuming the point values $\mU_j$ are available at a certain time level $t$, the solution is evolved to the next time level $t+\dt$
according to \eref{2.8} with the following numerical fluxes:
\begin{equation*}
\bm{{\cal F}}_\jph=\bm{{\cal L}}_\jph[\bm U(t)],\quad
\bm{{\cal F}}^{\rm I}_\jph=\bm{{\cal L}}_\jph\big[\bm U^{\rm I}\big],\quad
\bm{{\cal F}}^{\rm II}_\jph=\bm{{\cal L}}_\jph\big[\bm U^{\rm II}\big]-\bm\omega_\jph(t),
\end{equation*}
where
\begin{equation*}
\begin{aligned}
&\bm{{\cal L}}_\jph[\bm U]:=
\frac{1}{60}\big[\mF(\mU_{j+3})-8\mF(\mU_{j+2})+37\mF(\mU_{j+1})+37\mF(\mU_j)-8\mF(\mU_{j-1})+\mF(\mU_{j-2})\big],\\[0.5ex]
&\bm\omega_\jph(t):=
\frac{3\dx}{128\dt}\big[\mU_{j+3}(t)-5\mU_{j+2}(t)+10\mU_{j+1}(t)-10\mU_j(t)+5\mU_{j-1}(t)-\mU_{j-2}(t)\big].
\end{aligned}
\end{equation*}

In the next numerical example, we demonstrate that the described 1-D quasi-linear scheme achieves the fifth order of accuracy in capturing a
smooth solution provided $\dt\sim(\dx)^{\frac{5}{3}}$.

\paragraph{Example---1-D Accuracy Test.} We consider the following smooth initial data \cite{KKOKC}:
\begin{equation*}
u(x,0)=\sin\Big(\frac{\pi x}{5}+\frac{\pi}{4} \Big),\quad\rho(x,0)=
\bigg[\frac{\gamma-1}{2\sqrt{\gamma}}\left(u(x,0)+10\right)\bigg]^{\frac{2}{\gamma-1}},\quad p(x,0)=\rho^\gamma(x,0),
\end{equation*}
subject to the periodic boundary conditions in the computational domain $[0,10]$. We compute the numerical solution until the final time
$t=0.1$ on a sequence of uniform meshes with $\dx=1/10$, $1/20$, $1/40$, $1/80$, $1/160$, and $1/320$.

We then compute the $L^1$-errors and estimate the experimental convergence rates using the following Runge formulae, which are based on the
solutions computed on the three consecutive uniform grids with the mesh sizes $\dx$, $2\dx$, and $4\dx$ and denoted by $(\cdot)^{\dx}$,
$(\cdot)^{2\dx}$, and $(\cdot)^{4\dx}$, respectively:
$$
{\rm Error}(\dx)\approx\frac{\delta_{12}^2}{|\delta_{12}-\delta_{24}|},\quad
{\rm Rate}(\dx)\approx\log_2\left(\frac{\delta_{24}}{\delta_{12}}\right).
$$
Here, $\delta_{12}:=\|(\cdot)^{\dx}-(\cdot)^{2\dx}\|_{L^1}$ and $\delta_{24}:=\|(\cdot)^{2\dx}-(\cdot)^{4\dx}\|_{L^1}$. The obtained results
are reported in Table \ref{tab1}, where one can clearly see that the fifth order of accuracy is achieved. 
\begin{table}[ht!]
\centering
\begin{tabular}{|c|cc|cc|cc|cccc|cc|}
\hline
\multirow{2}{1em}{$\dx$}&\multicolumn{2}{c|}{$\rho$}&\multicolumn{2}{c|}{$\rho u$}&\multicolumn{2}{c|}{$E$}\\
\cline{2-7}&Error&Rate &Error&Rate&Error&Rate\\
\hline
$1/40$ &1.84e-07&4.76&6.03e-07&4.76&2.67e-06&4.76\\
$1/80$ &4.89e-09&4.99&1.59e-08&5.00&7.06e-08&5.00\\
$1/160$&1.56e-10&4.98&5.11e-10&4.98&2.26e-09&4.98\\
$1/320$&5.05e-12&4.97&1.65e-11&4.97&7.30e-11&4.97\\
\hline
\end{tabular}
\caption{\sf The $L^1$-errors and experimental convergence rates for $\rho$, $\rho u$, and $E$ computed by the quasi-linear fifth-order
scheme.\label{tab1}}
\end{table}

\subsection{1-D Adaptive Algorithm}\label{sec23}
We now turn to the description of the proposed 1-D adaptive scheme. To this end, we first compute the density-based SI proposed in
\cite{RL87} (see also \cite[\S2.2]{CKM2025}), which is based on the quantities
\begin{equation}
{\cal E}_j^\rho=\frac{|\rho_{j+1}-2\rho_j+\rho_{j-1}|}{|\rho_{j+1}-\rho_j|+|\rho_j-\rho_{j-1}|+\eps\big(\rho_{j+1}+2\rho_j+\rho_{j-1}\big)},
\label{2.9a}
\end{equation}
which we smooth out in the following way:
\begin{equation}
\xbar{\cal E}_j^{\,\rho}=\frac{1}{6}\big[{\cal E}_{j+1}^\rho+4{\cal E}_j^\rho+{\cal E}_{j-1}^\rho\big].
\label{2.9b}
\end{equation}
The term with $\eps$ in \eref{2.9a} plays the role of a ``noise'' filter, which is added in order not to refine ``wiggles'' or ``ripples''
that may appear due to loss of monotonicity; see \cite{RL87}. The value of $\eps$ depends on the problem at hand and in all of the numerical
examples presented in this paper, we use the same value $\eps=0.2$ as in \cite{RL87,CKM2025}. We then compute the
pressure-based SI, which is calculated using the following quantities, which are similar to \eref{2.9a}--\eref{2.9b}:
$$
{\cal E}_j^p=\frac{|p_{j+1}-2p_j+p_{j-1}|}{|p_{j+1}-p_j|+|p_j-p_{j-1}|+\eps\big(p_{j+1}+2p_j+p_{j-1}\big)},\quad
\xbar{\cal E}_j^{\,p}=\frac{1}{6}\big[{\cal E}_{j+1}^p+4{\cal E}_j^p+{\cal E}_{j-1}^p\big].
$$

As in \cite{CKM2025}, we identify the ``rough'' parts of the computed solution with the help of $\,\xbar{\cal E}^{\,\rho}_j$,
which are large wherever the solution is nonsmooth. To distinguish between shock and contact waves, we take advantage of the fact that, at
contact discontinuities, $\,\xbar{\cal E}^{\,\rho}_j$ is large, but $\,\xbar{\cal E}^{\,p}_j$ is small, and design the following adaptive
strategy:

\smallskip
\noindent
$\bullet$ Select a tunable constant $0<\texttt{C}_1<1$ and mark the cells $I_j$, in which $\,\xbar{\cal E}^{\,\rho}_j>\texttt{C}_1$, as
``rough'' cells;

\smallskip
\noindent
$\bullet$ Select a constant $0<\texttt{C}_2<1$ and mark those ``rough'' cells $I_j$, in which $\,\xbar{\cal E}^{\,p}_j<\texttt{C}_2$, as
``contact'' cells.

After identifying three different parts of the computed solution, we design a novel scheme adaption approach as follows:

\smallskip
\noindent
\underline{Area A} (neighborhoods of the contact discontinuities): Use the second-order LDCU scheme (\S\ref{sec21}) with the overcompressive
SBM limiter (Appendix \ref{sec211}) with $\tau=-0.25$ in \eref{2.6};

\smallskip
\noindent
\underline{Area B} (the rest of the ``rough'' parts of the computed solution): Use the second-order LDCU scheme (\S\ref{sec21}) with the
Minmod2 limiter, which is, in fact, a dissipative SBM limiter (Appendix \ref{sec211}) with $\tau=0.5$ in \eref{2.6};

\smallskip
\noindent
\underline{Area C} (smooth parts of the computed solution): Use the quasi-linear fifth-order scheme (\S\ref{sec22}).

\section{Novel 2-D Adaptive Scheme}\label{sec3}
In this section, we extend the novel 1-D adaptive strategy introduced in \S\ref{sec2} to the 2-D case.

\subsection{2-D Low-Dissipation Central-Upwind (LDCU) Schemes}\label{sec31}
In this section, we briefly review the second-order semi-discrete 2-D LDCU schemes introduced in \cite{CKX_24}.

We assume that the computational domain is covered with uniform cells $I_{j,k}:=\big[x_\jmh,x_\jph\big]\times\big[y_\kmh,y_\kph\big]$ with
$x_\jph-x_\jmh\equiv\dx$ and $y_\kph-y_\kmh\equiv\dy$ centered at $(x_j,y_k)$ with $x_j=\big(x_\jmh+x_\jph\big)/2$ and
$y_k=\big(y_\kmh+y_\kph\big)/2$, and that the cell averages
\begin{equation*}
\xbar\mU_{j,k}:\approx\frac{1}{\dx\dy}\iint\limits_{I_{j,k}}\mU(x,y,t)\,{\rm d}x\,{\rm d}y
\end{equation*}
are available at a certain time $t\ge0$.

According to \cite{CKX_24}, $\,\xbar\mU_{j,k}$ are evolved in time by numerically solving the following system of ODEs:
\begin{equation}
\frac{{\rm d}\xbar\mU_{j,k}}{{\rm d}t}=-\frac{\bmF^{\rm LDCU}_{\jph,k}[\bm U]-\bmF^{\rm LDCU}_{\jmh,k}[\bm U]}{\dx}-
\frac{\bmG^{\rm LDCU}_{j,\kph}[\bm U]-\bmG^{\rm LDCU}_{j,\kmh}[\bm U]}{\dy},
\label{4.1}
\end{equation}
where $\bm{{\cal F}}^{\rm LDCU}_{\jph,k}[\bm U]$ and $\bm{{\cal G}}^{\rm LDCU}_{j,\kph}[\bm U]$ are the $x$- and $y$-directional LDCU numerical fluxes given by
\begin{equation*}
\begin{aligned}
\bmF^{\rm LDCU}_{\jph,k}[\bm U]&=\frac{a^+_{\jph,k}\mF\big(\bm U^-_{\jph,k}\big)-a^-_{\jph,k}\mF\big(\bm U^+_{\jph,k}\big)}
{a^+_{\jph,k}-a^-_{\jph,k}}+\frac{a^+_{\jph,k}a^-_{\jph,k}}{a^+_{\jph,k}-a^-_{\jph,k}}\Big(\mU^+_{\jph,k}-\mU^-_{\jph,k}\Big)+
\bm q^x_{\jph,k},\\
\bmG^{\rm LDCU}_{j,\kph}[\bm U]&=\frac{b^+_{j,\kph}\mG\big(\bm U^-_{j,\kph}\big)-b^-_{j,\kph}\mG\big(\bm U^+_{j,\kph}\big)}
{b^+_{j,\kph}-b^-_{j,\kph}}+\frac{b^+_{j,\kph}b^-_{j,\kph}}{b^+_{j,\kph}-b^-_{j,\kph}}\Big(\mU^+_{j,\kph}-\mU^-_{j,\kph}\Big)+
\bm q^y_{j,\kph}.
\end{aligned}
\end{equation*}
Here:

\smallskip
\noindent
$\bullet$ $\mU^\pm_{\jph,k}$ and $\mU^\pm_{j,\kph}$ are the one-sided point values of $\mU$ at the midpoints of the cell interfaces
$(x_\jph,y_k)$ and $(x_j,y_\kph)$, respectively. As in the 1-D case, these one-sided point values are reconstructed out of the given set of
cell averages $\{\,\xbar\mU_{j,k}\}$ using a proper nonlinear limiter; see Appendix \ref{appa};

\smallskip
\noindent
$\bullet$ $a^\pm_{\jph,k}$ and $b^\pm_{j,\kph}$ are the one-sided local speeds of propagation in the $x$- and $y$-directions, respectively.
They can be estimated using the largest and smallest eigenvalues of the corresponding Jacobians
$A(\mU):=\frac{\partial\mF}{\partial\mU}(\mU)$ and $B(\mU):=\frac{\partial\mG}{\partial\mU}(\mU)$;

\smallskip
\noindent
$\bullet$ $\bm q^x_{\jph,k}$ and $\bm q^y_{j,\kph}$ are ``built-in'' anti-diffusion terms, which need to be derived for a particular system
\eref{1.2} at hand.

\smallskip
As in the 1-D case, we focus on the Euler equations of gas dynamics, whose 2-D version reads as \eref{1.2} with
\begin{equation}
\mU=\big(\rho,\rho u,\rho v,E\big)^\top,\quad\mF=\big(\rho u,\rho u^2+p,\rho uv,u(E+p)\big)^\top,\quad
\mG=\big(\rho v,\rho uv,\rho v^2+p,v(E+p)\big)^\top.
\label{3.3}
\end{equation}
Here, $v$ is the $y$-velocity and the other variables are the same as in the 1-D case. The system \eref{1.2}, \eref{3.3} is closed through
the following EOS for ideal gases:
\begin{equation}
p=(\gamma-1)\Big[E-\frac{\rho}{2}(u^2+v^2)\Big].
\label{3.4}
\end{equation}
For the Euler system \eref{1.2}, \eref{3.3}--\eref{3.4}, the one-sided local speeds of propagation can be estimated by
\begin{equation*}\resizebox{\linewidth}{!}{$
\begin{aligned}
&a^+_{\jph,k}=\max\Big\{u^+_{\jph,k}+c^+_{\jph,k},u^-_{\jph,k}+c^-_{\jph,k},0\Big\},~
&a^-_{\jph,k}=\min\Big\{u^+_{\jph,k}-c^+_{\jph,k},u^-_{\jph,k}-c^-_{\jph,k},0\Big\},\\
&b^+_{j,\kph}=\max\Big\{v^+_{j,\kph}+c^+_{j,\kph},v^-_{j,\kph}+c^-_{j,\kph},0\Big\},~
&b^-_{j,\kph}=\min\Big\{v^+_{j,\kph}-c^+_{j,\kph},v^-_{j,\kph}-c^-_{j,\kph},0\Big\},\\
\end{aligned}$}
\end{equation*}
where
$$
\begin{aligned}
&u^\pm_{\jph,k}=\frac{(\rho u)_{\jph,k}^\pm}{\rho^\pm_{\jph,k}},\quad u^\pm_{j,\kph}=\frac{(\rho u)_{j,\kph}^\pm}{\rho^\pm_{j,\kph}},\quad
v^\pm_{\jph,k}=\frac{(\rho v)_{\jph,k}^\pm}{\rho^\pm_{\jph,k}},\quad v^\pm_{j,\kph}=\frac{(\rho v)_{j,\kph}^\pm}{\rho^\pm_{j,\kph}},\\
&p^\pm_{\jph,k}=(\gamma-1)\Big[E^\pm_{\jph,k}-\hf\rho^\pm_{\jph,k}\Big(\big(u^\pm_{\jph,k}\big)^2+\big(v^\pm_{\jph,k}\big)^2\Big)\Big],\quad
c^\pm_{\jph,k}=\sqrt{\frac{\gamma p^\pm_{\jph,k}}{\rho^\pm_{\jph,k}}},\\
&p^\pm_{j,\kph}=(\gamma-1)\Big[E^\pm_{j,\kph}-\hf\rho^\pm_{j,\kph}\Big(\big(u^\pm_{j,\kph}\big)^2+\big(v^\pm_{j,\kph}\big)^2\Big)\Big],\quad
c^\pm_{j,\kph}=\sqrt{\frac{\gamma p^\pm_{j,\kph}}{\rho^\pm_{j,\kph}}}.
\end{aligned}
$$
Finally, a description of the anti-diffusion terms $\bm q^x_{\jph,k}$ and $\bm q^y_{j,\kph}$ can be found in \cite{CKX_24}, and here, we
omit the details for the sake of brevity.

As in the 1-D case, we numerically integrate the ODE system \eref{4.1} using the three-stage third-order SSP Runge-Kutta method:
\begin{equation}
\begin{aligned}
&\mU^{\rm I}_{j,k}=\mU_{j,k}(t)-\frac{\dt}{\dx}\left(\bm{{\cal F}}_{\jph,k}-\bm{{\cal F}}_{\jmh,k}\right)
-\frac{\dt}{\dy}\left(\bm{{\cal G}}_{j,\kph}-\bm{{\cal G}}_{j,\kmh}\right),\\
&\mU^{\rm II}_{j,k}=\frac{3}{4}\,\mU_{j,k}(t)+\frac{1}{4}\left[\mU^{\rm I}_{j,k}-
\frac{\dt}{\dx}\left(\bm{{\cal F}}^{\rm I}_{\jph,k}-\bm{{\cal F}}^{\rm I}_{\jmh,k}\right)\right.\left.-\frac{\dt}{\dy}\left(\bm{{\cal G}}^{\rm I}_{j,\kph}-
\bm{{\cal G}}^{\rm I}_{j,\kmh}\right)\right],\\
&\mU_{j,k}(t+\dt)=\frac{1}{3}\,\mU_{j,k}(t)+\frac{2}{3}\left[\mU^{\rm II}_{j,k}-
\frac{\dt}{\dx}\left(\bm{{\cal F}}^{\rm II}_{\jph,k}-\bm{{\cal F}}^{\rm II}_{\jmh,k}\right)\right.
\left.-\frac{\dt}{\dy}\left(\bm{{\cal G}}^{\rm II}_{j,\kph}-\bm{{\cal G}}^{\rm II}_{j,\kmh}\right)\right].
\end{aligned}
\label{4.4}
\end{equation}
where $\bm{{\cal F}}_{\jph,k}=\bm{{\cal F}}^{\rm LDCU}_{\jph,k}[\mU]$, $\bm{{\cal F}}^{\rm I}_{\jph,k}=\bm{{\cal F}}^{\rm LDCU}_{\jph,k}[\mU^{\rm I}]$, $\bm{{\cal F}}^{\rm II}_{\jph,k}=\bm{{\cal F}}^{\rm LDCU}_{\jph,k}[\mU^{\rm II}]$, $\bm{{\cal G}}_{j,\kph}=\bm{{\cal G}}^{\rm LDCU}_{j,\kph}[\mU]$, $\bm{{\cal G}}^{\rm I}_{j,\kph}=\bm{{\cal G}}^{\rm LDCU}_{j,\kph}[\mU^{\rm I}]$, and $\bm{{\cal G}}^{\rm II}_{j,\kph}=\bm{{\cal G}}^{\rm LDCU}_{j,\kph}[\mU^{\rm II}]$.

In \eref{4.4}, the time step $\dt$ is selected adaptively based on the following CFL-based stability restriction:
\begin{equation*}
\dt\le\min\left\{\frac{\dx}{2a},\frac{\dy}{2b}\right\},\quad a:=\max_{j,k}\Big\{\max\big(a^+_{\jph,k},-a^-_{\jph,k}\big)\Big\},\quad
b:=\max_{j,k}\Big\{\max\big(b^+_{j,\kph},-b^-_{j,\kph}\big)\Big\}.
\end{equation*}

As in the 1-D case, we have replaced in \eref{4.4} the cell averages of $\bm U$ with its point values
$\bm U_{j,k}(t):\approx\bm U(x_j,y_k,t)$. This does not affect the accuracy of the second-order scheme as the difference between the
cell averages and point values is proportional to ${\cal O}((\dx)^2+(\dy)^2)$ for smooth solutions.

\subsection{2-D Quasi-Linear Fifth-Order Scheme}\label{sec42}
In this section, we extend the 1-D quasi-linear fifth-order scheme described in \S\ref{sec22} to the 2-D case. The extension is performed in
a ``dimension-by-dimension'' manner.

Assuming the point values $\mU_{j,k}$ are available at a certain time level $t$, the solution is evolved to the next time level $t+\dt$
according to \eref{4.4} with the following numerical fluxes:
\begin{equation*}
\begin{aligned}
&\bm{{\cal F}}_{\jph,k}=\bm{{\cal L}}_{\jph,k}^x[\bm U(t)],&&
\bm{{\cal F}}_{\jph,k}^{\rm I}=\bm{{\cal L}}_{\jph,k}^x\big[\bm U^{\rm I}\big],&&
\bm{{\cal F}}_{\jph,k}^{\rm II}=\bm{{\cal L}}_{\jph,k}^x\big[\bm U^{\rm II}\big]-\bm\omega^x_{\jph,k}(t),\\
&\bm{{\cal G}}_{j,\kph}=\bm{{\cal L}}_{j,\kph}^y[\bm U(t)],&&
\bm{{\cal G}}_{j,\kph}^{\rm I}=\bm{{\cal L}}_{j,\kph}^y\big[\bm U^{\rm I}\big],&&
\bm{{\cal G}}_{j,\kph}^{\rm II}=\bm{{\cal L}}_{j,\kph}^y\big[\bm U^{\rm II}\big]-\bm\omega^y_{j,\kph}(t),
\end{aligned}
\end{equation*}
where
\begin{equation*}\resizebox{\linewidth}{!}{$
\begin{aligned}
&\bm{{\cal L}}_{\jph,k}^x[\bm U]:=
\frac{1}{60}\big[\mF(\mU_{j+3,k})-8\mF(\mU_{j+2,k})+37\mF(\mU_{j+1,k})+37\mF(\mU_{j,k})-8\mF(\mU_{j-1,k})+\mF(\mU_{j-2,k})\big],\\[0.5ex]
&\bm\omega_{\jph,k}^x(t):=
\frac{3\dx}{128\dt}\big[\mU_{j+3,k}(t)-5\mU_{j+2,k}(t)+10\mU_{j+1,k}(t)-10\mU_{j,k}(t)+5\mU_{j-1,k}(t)-\mU_{j-2,k}(t)\big],\\[0.5ex]
&\bm{{\cal L}}_{j,\kph}^y[\bm U]:=
\frac{1}{60}\big[\mG(\mU_{j,k+3})-8\mG(\mU_{j,k+2})+37\mG(\mU_{j,k+1})+37\mG(\mU_{j,k})-8\mG(\mU_{j,k-1})+\mG(\mU_{j,k-2})\big],\\[0.5ex]
&\bm\omega_{j,\kph}^y(t):=
\frac{3\dy}{128\dt}\big[\mU_{j,k+3}(t)-5\mU_{j,k+2}(t)+10\mU_{j,k+1}(t)-10\mU_{j,k}(t)+5\mU_{j,k-1}(t)-\mU_{j,k-2}(t)\big],
\end{aligned}$}
\end{equation*}

In the next numerical example, we demonstrate that the developed 2-D quasi-linear scheme achieves the fifth order of accuracy in capturing a
smooth solution.

\paragraph{Example---2-D Accuracy Test.} We consider the following smooth initial data \cite{BD2013,Shu1998}:
\begin{equation*}
\begin{aligned}
&\rho(x,y,0)=\bigg(1-\frac{(\gamma-1)\kappa^2}{2\gamma}\bigg)^{\frac{1}{\gamma-1}},\quad p(x,y,0)=\rho^\gamma(x,y,0),\\
&u(x,y,0)=1-\kappa y, \quad v(x,y,0)=1+\kappa x, \quad \kappa = \frac{5}{2 \pi} e^{\frac{1-x^2-y^2}{2}} 
\end{aligned}
\end{equation*}
subject to the periodic boundary conditions in the computational domain $[-10,10]\times[-10,10]$. The exact solution of this initial value
problem is given by $\mU(x,y,t)=\mU(x-t,y-t,0)$.

We compute the numerical solution until the final time $t=0.1$ on a sequence of uniform meshes with $\dx=\dy=1/10$, $1/20$, $1/40$,
and $1/80$ and $\dt\sim\min\{(\dx)^{\frac{5}{3}},(\dy)^{\frac{5}{3}} \}$, and then measure the $L^1$-errors and the corresponding
experimental convergence rates. The obtained results are presented in Table \ref{tab52}, where one can see that the expected fifth order of
accuracy has been achieved.
\begin{table}[ht!]
\centering
\begin{tabular}{|c|cc|cc|cc|cc|cc|cc|}
\hline
\multirow{2}{5em}{$\dx=\dy$}&\multicolumn{2}{c|}{$\rho$}&\multicolumn{2}{c|}{$u$}&\multicolumn{2}{c|}{$v$}&\multicolumn{2}{c|}{$p$}\\
\cline{2-9}&Error&Rate &Error&Rate&Error&Rate &Error&Rate\\
\hline
$1/10$ &3.35e-04&--- &6.01e-04&--- &6.33e-04&--- &4.85e-04&---\\
$1/20$ &1.02e-05&5.03&1.82e-05&5.04&1.97e-05&5.01&1.48e-05&5.04\\
$1/40$ &3.56e-07&4.85&6.36e-07&4.84&6.89e-07&4.84&5.14e-07&4.85\\
$1/80$ &1.15e-08&4.95&2.07e-08&4.94&2.24e-08&4.94&1.66e-08&4.95\\
\hline
\end{tabular}
\caption{\sf The $L^1$-errors and experimental convergence rates for $\rho$, $u$, $v$, and $p$ computed by the quasi-linear fifth-order
scheme.\label{tab52}}
\end{table}

\subsection{2-D Adaptive Algorithm}
We now turn to the description of the proposed 2-D adaptive schemes.

We first introduce the quantities \cite{CKM2025}, which are similar to \eref{2.9a}:
\begin{equation*}
\begin{aligned}
{\cal E}^\rho_{j,k}=\sqrt{\frac{{\cal E}^{1,\rho}_{j,k}}{{\cal E}^{2,\rho}_{j,k}}},\quad
{\cal E}^{1,\rho}_{j,k}&=\big(\rho_{j+1,k}-2\rho_{j,k}+\rho_{j-1,k}\big)^2+\big(\rho_{j,k+1}-2\rho_{j,k}+\rho_{j,k-1}\big)^2,\\
{\cal E}^{2,\rho}_{j,k}&=\big(|\rho_{j+1,k}-\rho_{j,k}|+|\rho_{j,k}-\rho_{j-1,k}|+
\varepsilon\big[|\xbar\rho_{j+1,k}|+2|\rho_{j,k}|+|\rho_{j-1,k}|\big]\big)^2\\
&+\big(|\rho_{j,k+1}-\rho_{j,k}|+|\rho_{j,k}-\rho_{j,k-1}|+\varepsilon\big[|\rho_{j,k+1}|+2|\rho_{j,k}|+|\rho_{j,k-1}|\big]\big)^2,
\end{aligned}
\end{equation*}
and to \eref{2.9b}:
\begin{equation*}
\xbar{\cal E}^{\,\rho}_{j,k}=\frac{1}{36}\Big[{\cal E}^{\,\rho}_{j-1,k-1}+{\cal E}^{\,\rho}_{j-1,k+1}+{\cal E}^{\,\rho}_{j+1,k-1}+
{\cal E}^{\,\rho}_{j+1,k+1}+4\big({\cal E}^{\,\rho}_{j-1,k}+{\cal E}^{\,\rho}_{j,k-1}+{\cal E}^{\,\rho}_{j,k+1}+
{\cal E}^{\,\rho}_{j+1,k}\big)+16{\cal E}^{\,\rho}_{j,k}\Big].
\end{equation*}
The values of ${\cal E}^p_{j,k}$ and $\,\xbar{\cal E}^{\,p}_{j,k}$ can be computed in a same manner by replacing $\rho$ with $p$.

Equipped with $\,\xbar{\cal E}^{\,\rho}_{j,k}$ and $\,\xbar{\cal E}^{\,p}_{j,k}$, we proceed as in the 1-D case:

As in the 1-D case, we use an overcompressive SBM
limiter in the contact wave regions of the computed solution, the Minmod2 limiter in the shock wave regions, and the fifth-order scheme in the rest parts. To this end, we detect the smoothness of the numerical solution using the 2-D version of the SI proposed in \cite{RL87}:

\smallskip
\noindent
$\bullet$ Select a tunable constant $0<\texttt{C}_1<1$ and mark the cells $I_{j,k}$, in which $\,\xbar{\cal E}^{\,\rho}_{j,k}>\texttt{C}_1$,
as ``rough'' cells;

\smallskip
\noindent
$\bullet$ Select a constant $0<\texttt{C}_2<1$ and mark those ``rough'' cells $I_{j,k}$, in which
$\,\xbar{\cal E}^{\,p}_{j,k}<\texttt{C}_2$, as ``contact'' cells.

After identifying three different parts of the computed solution, we design a 2-D scheme adaption algorithm as follows:

\smallskip
\noindent
\underline{Area A} (neighborhoods of the contact discontinuities): Use the second-order LDCU scheme (\S\ref{sec31}) with the overcompressive
SBM limiter (Appendix \ref{appa}) with $\tau=-0.25$ in \eref{2.6};

\smallskip
\noindent
\underline{Area B} (the rest of the ``rough'' parts of the computed solution): Use the second-order LDCU scheme (\S\ref{sec31}) with the
Minmod2 limiter, which is, in fact, a dissipative SBM limiter (Appendix \ref{appa}) with $\tau=0.5$ in \eref{2.6};

\smallskip
\noindent
\underline{Area C} (smooth parts of the computed solution): Use the quasi-linear fifth-order scheme (\S\ref{sec42}).

\section{Numerical Examples}\label{sec4}
In this section, we apply the developed adaptive schemes to both the 1-D and 2-D Euler equations of gas dynamics and compare the performance
of the adaptive schemes from \cite{CKM2025} (OLD schemes) and novel adaptive schemes, which will be referred to as NEW schemes. 
We take $\gamma=1.4$ in Examples 1--8 and $\gamma=5/3$ in Example 9. In all of the examples, we use the CFL number $0.4$.

All of the numerical experiments have been conducted on a workstation equipped with an Intel(R) Core(TM) i7-9750H CPU at 2.6 GHz and 32 GB
of RAM. The simulations were conducted in FORTRAN using GCC version 14.2.0 compiler suite. The reported CPU times were averaged over thirty
independent runs to ensure reproducibility and minimize variability due to system processes.

\subsection{1-D Examples}
In the 1-D examples, the adaptive solutions will be plotted along with the reference solutions, which will be computed using the
second-order LDCU scheme described in \S\ref{sec21} and implemented with the Minmod2 limiter.

\paragraph{Example 1---Shock-Density Wave Interaction Problem.} In the first example taken from \cite{SO89}, we consider the
shock-density wave interaction problem. The initial data,
\begin{equation*}
(\rho,u,p)\Big|_{(x,0)}=\begin{cases}\bigg(\dfrac{27}{7},\dfrac{4\sqrt{35}}{9},\dfrac{31}{3}\bigg),&x<-4,\\[0.8ex]
(1+0.2\sin(5x),0,1),&x>-4,
\end{cases}
\end{equation*}
are prescribed in the computational domain $[-5,15]$ subject to the free boundary conditions.

\begin{figure}[ht!]
\centerline{\includegraphics[trim=0.9cm 0.4cm 0.9cm 0.6cm, clip, width=6.0cm]{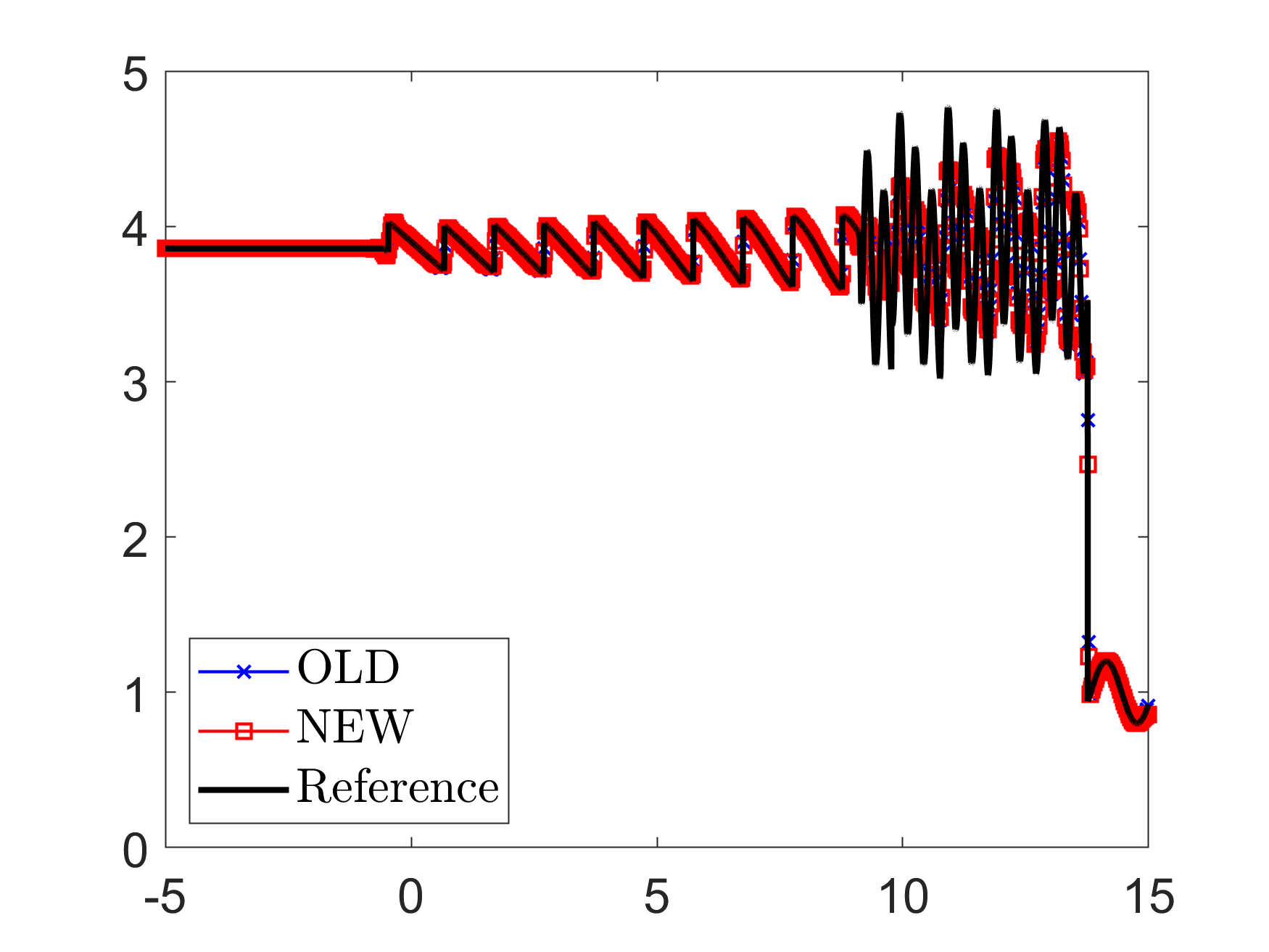}\hspace{1cm}
            \includegraphics[trim=0.9cm 0.4cm 0.9cm 0.6cm, clip, width=6.0cm]{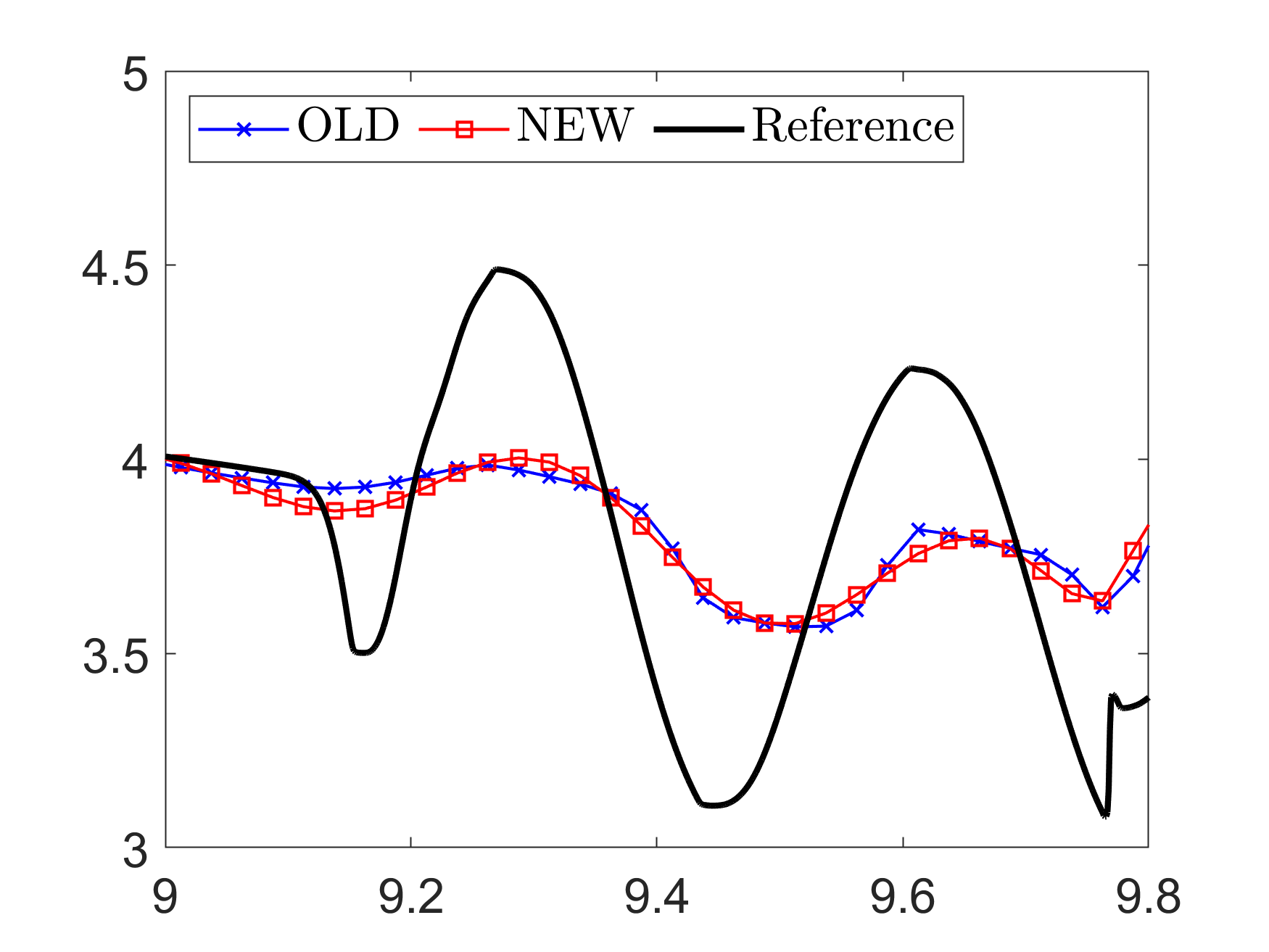}}
\caption{\sf Example 1: Density $\rho$ computed on a uniform mesh with $\dx=1/40$ by the NEW and OLD schemes (left) and zoom at
$x\in[9,9.8]$ (right).\label{fig4}}
\end{figure}

\begin{figure}[ht!]
\centerline{\includegraphics[trim=0.9cm 0.4cm 0.9cm 0.6cm, clip, width=6.0cm]{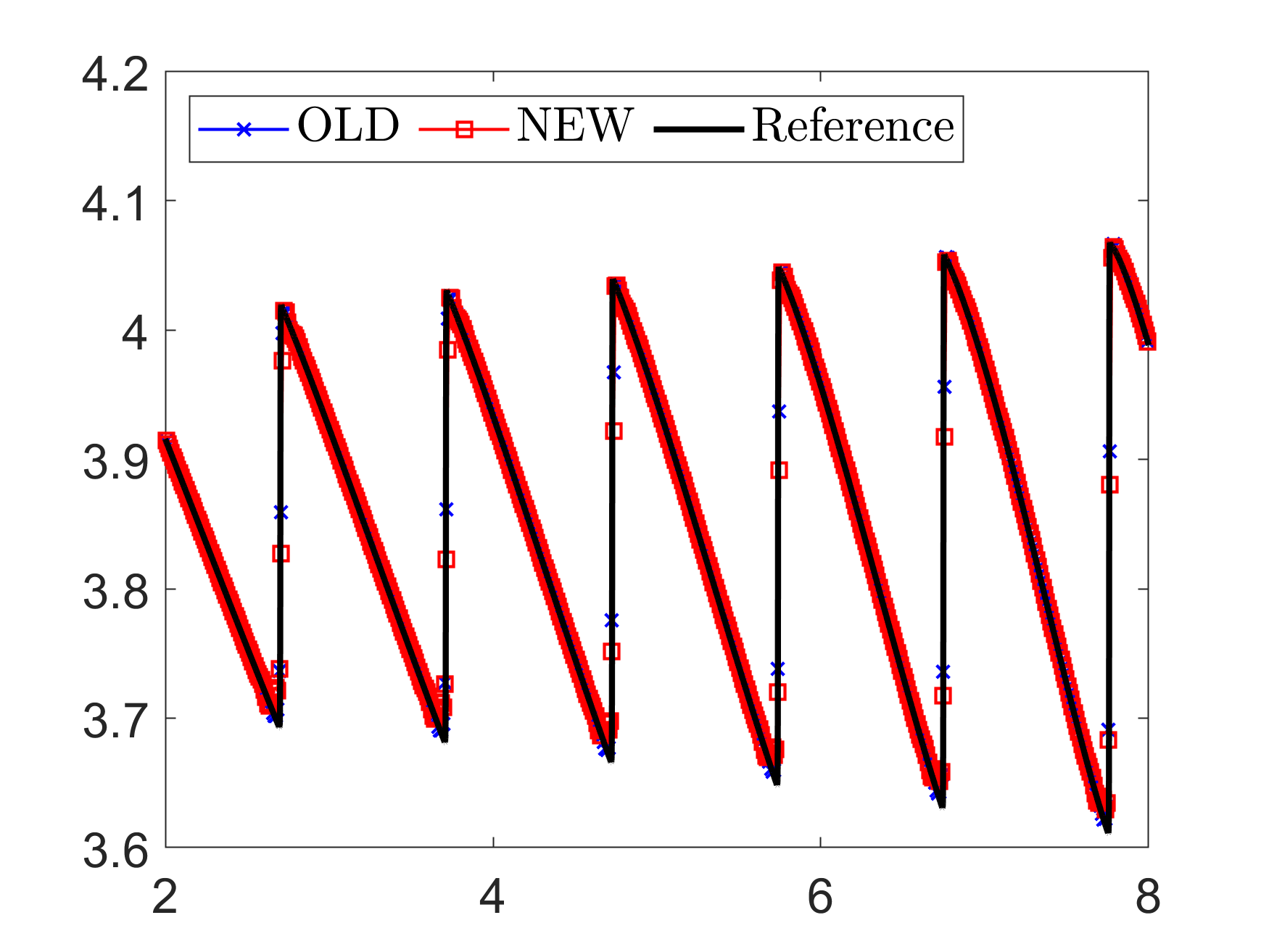}\hspace{1cm}
            \includegraphics[trim=0.9cm 0.4cm 0.9cm 0.6cm, clip, width=6.0cm]{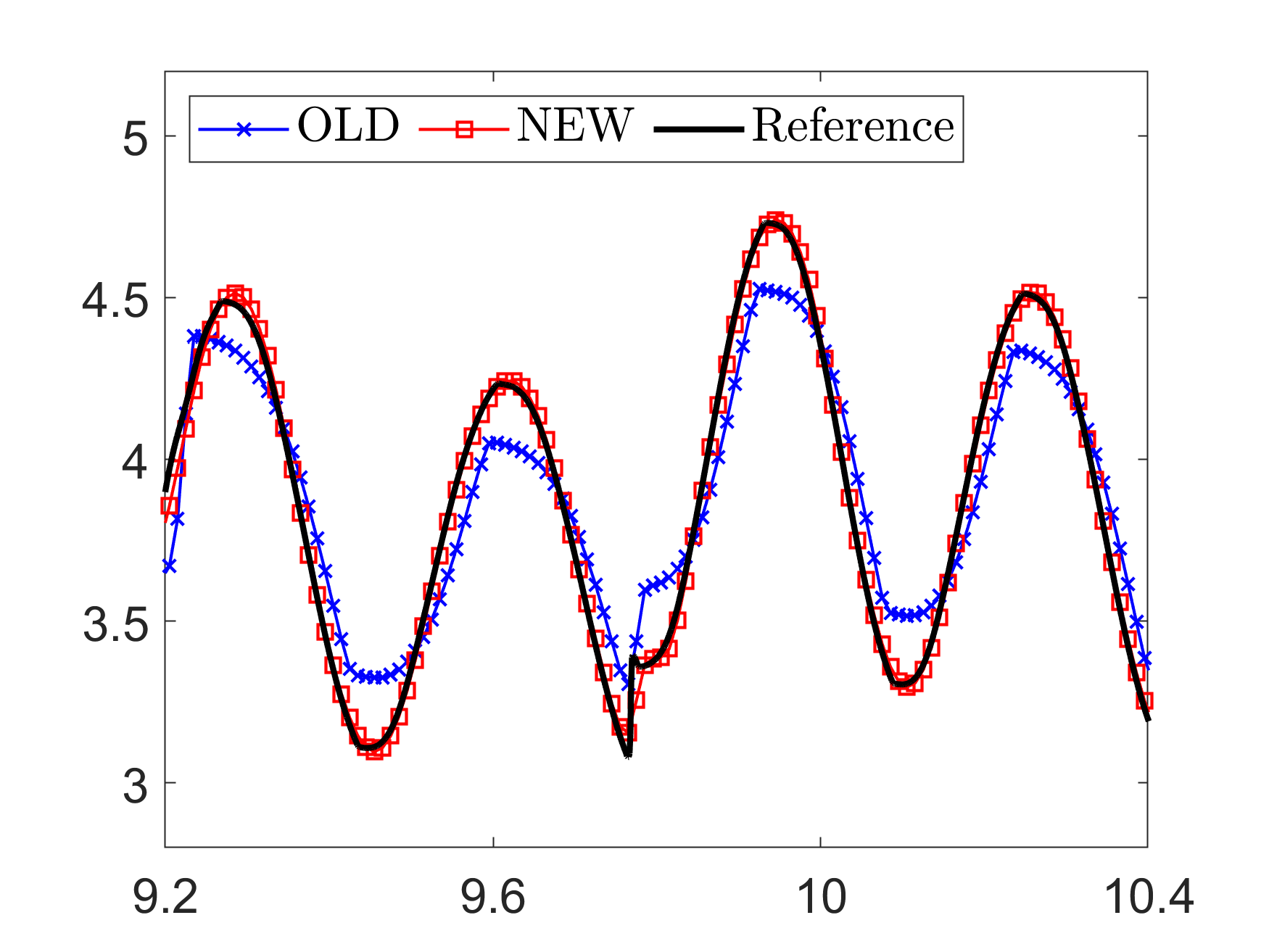}}
\caption{\sf Example 1: Density $\rho$ computed on a uniform mesh with $\dx=1/100$ by the NEW and OLD schemes---zoom at $x\in[4,8]$ (left)
and $x\in[9.2,10.4]$ (right).\label{fig4a}}
\end{figure}

We compute the numerical solutions by the NEW (with the adaption constants $\texttt{C}_1=0.015$ and $\texttt{C}_2=0.15$) and OLD (with the
adaption constant $\texttt{C}_1=0.01$) schemes on uniform meshes with $\dx=1/40$ and $1/100$ until the final time $t=5$. We present the
obtained numerical results in Figure \ref{fig4} together with the reference solution computed on a much finer mesh with $\dx=1/1000$. One
can see that, compared with the OLD scheme, the NEW one produces a slightly more accurate result; see Figure \ref{fig4} (right). At the same
time, there are neither oscillations near the shocks (see Figure \ref{fig4a} (left)) nor ``stair-like'' structures in the smooth parts of
the computed solution (see Figure \ref{fig4a} (right)), when a finer mesh with $\dx=1/200$ is used.

We then measure the CPU times consumed by both of the studied schemes on a uniform mesh with $\dx=1/4000$ and the CPU time consumed by the
OLD scheme is about $356\%$ larger than the CPU time consumed by the NEW scheme. This demonstrates that the NEW scheme is not only more
accurate, but also substantially more efficient than its OLD counterpart. We stress that the NEW scheme is highly efficient due to an
extremely low computational cost of the quasi-linear fifth-order scheme.

As mentioned in the Introduction, a simpler adaptive strategy, in which we simply replace in the smooth regions the LDCU scheme equipped
with the dissipative Minmod2 limiter (as in the OLD scheme from \cite{CKM2025}) by the quasi-linear fifth-order scheme from \cite{KOK25},
may fail to produce satisfactory results. In such a simplified version of the modification of the OLD scheme, which will be referred to as
the MODIFIED scheme, we use the same way to detect ``rough'' parts of the computed solutions as in \cite{CKM2025}: If
$\,\xbar{\cal E}_j^{\,\rho}\ge\texttt{C}_1$, the cell $I_j$ is marked as ``rough''.

\begin{figure}[ht!]
\centerline{\includegraphics[trim=0.9cm 0.4cm 0.8cm 0.6cm, clip, width=6.0cm]{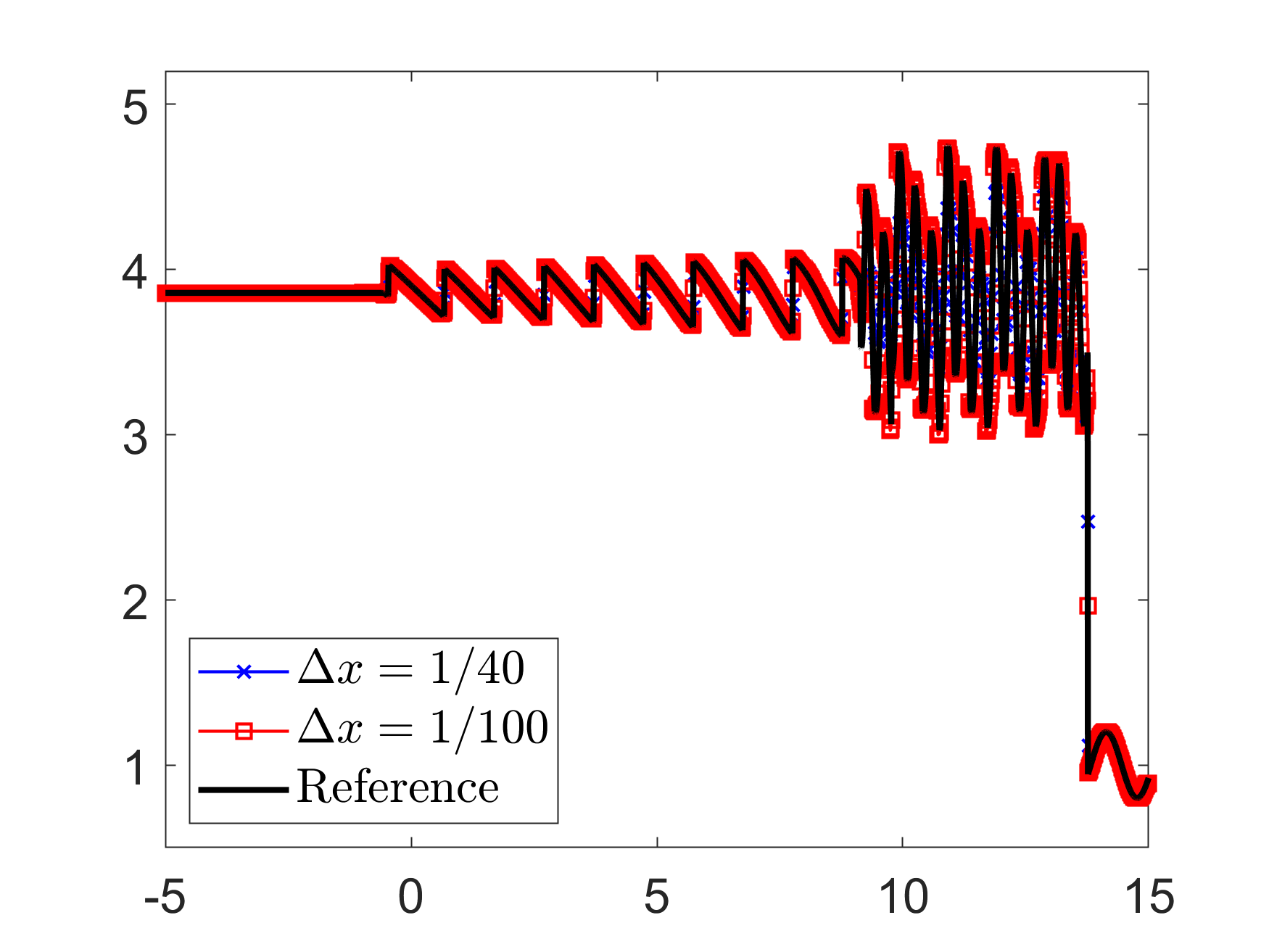}\hspace{1cm}
            \includegraphics[trim=0.9cm 0.4cm 0.8cm 0.6cm, clip, width=6.0cm]{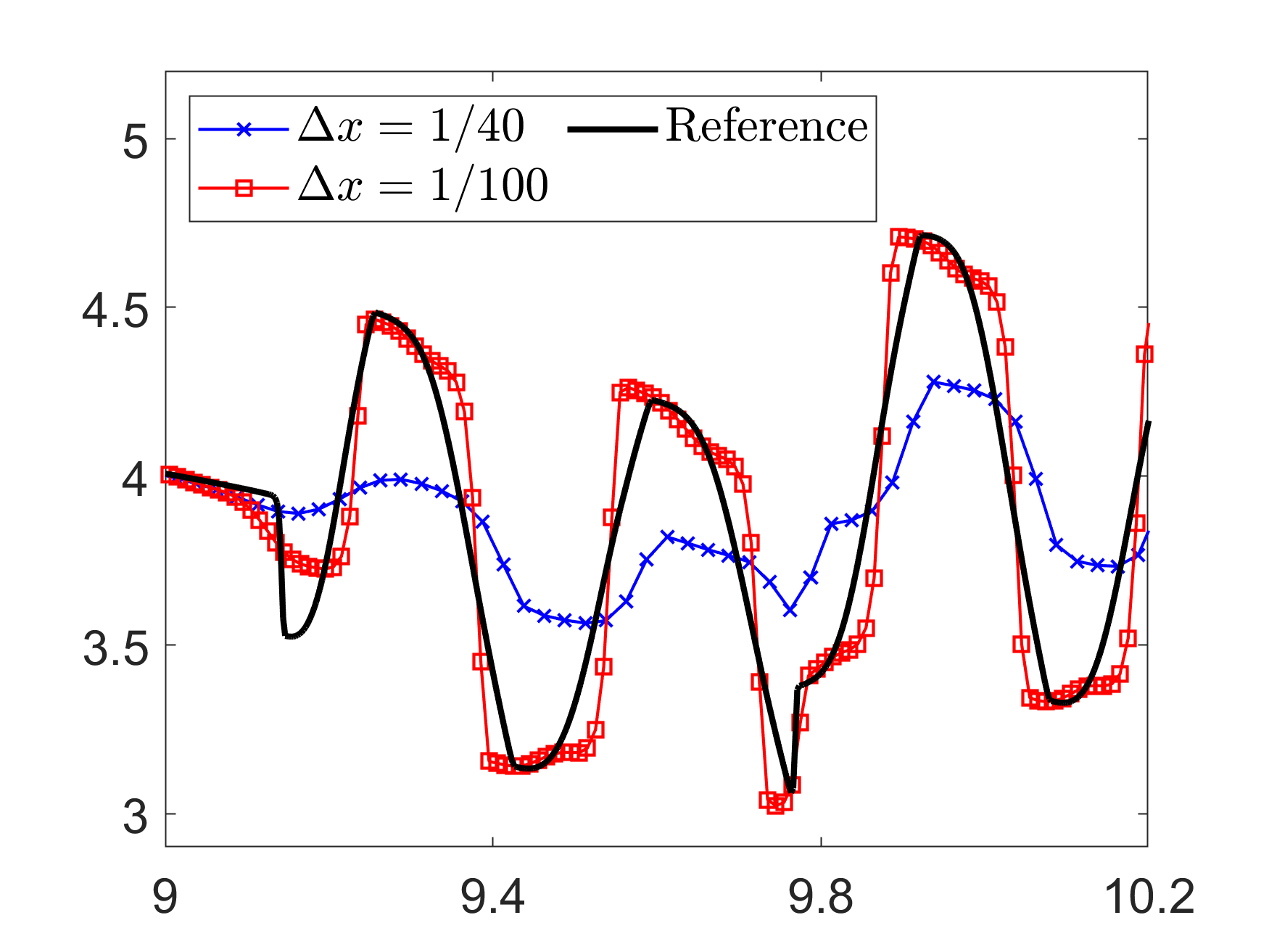}}
\caption{\sf Example 1: Density $\rho$ computed with $\dx=1/40$ and $\dx=1/100$ by the MODIFIED scheme (left) and zoom at
$x\in[9,10.2]$ (right).\label{fig1a}}
\end{figure}

\begin{figure}[ht!]
\centerline{\includegraphics[trim=0.9cm 0.4cm 1.1cm 0.6cm, clip, width=6.0cm]{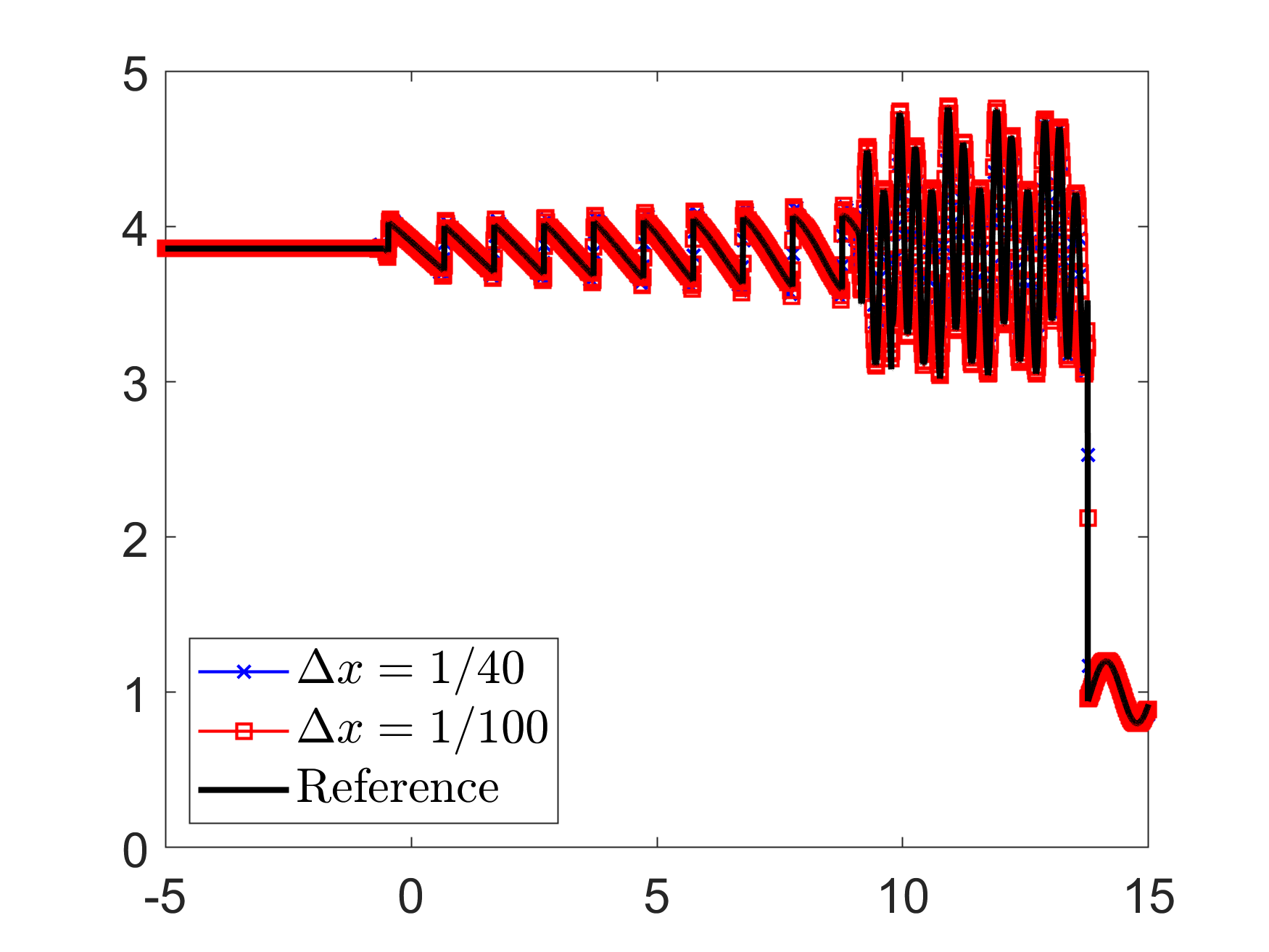}\hspace{1cm}
            \includegraphics[trim=0.9cm 0.4cm 1.1cm 0.6cm, clip, width=6.0cm]{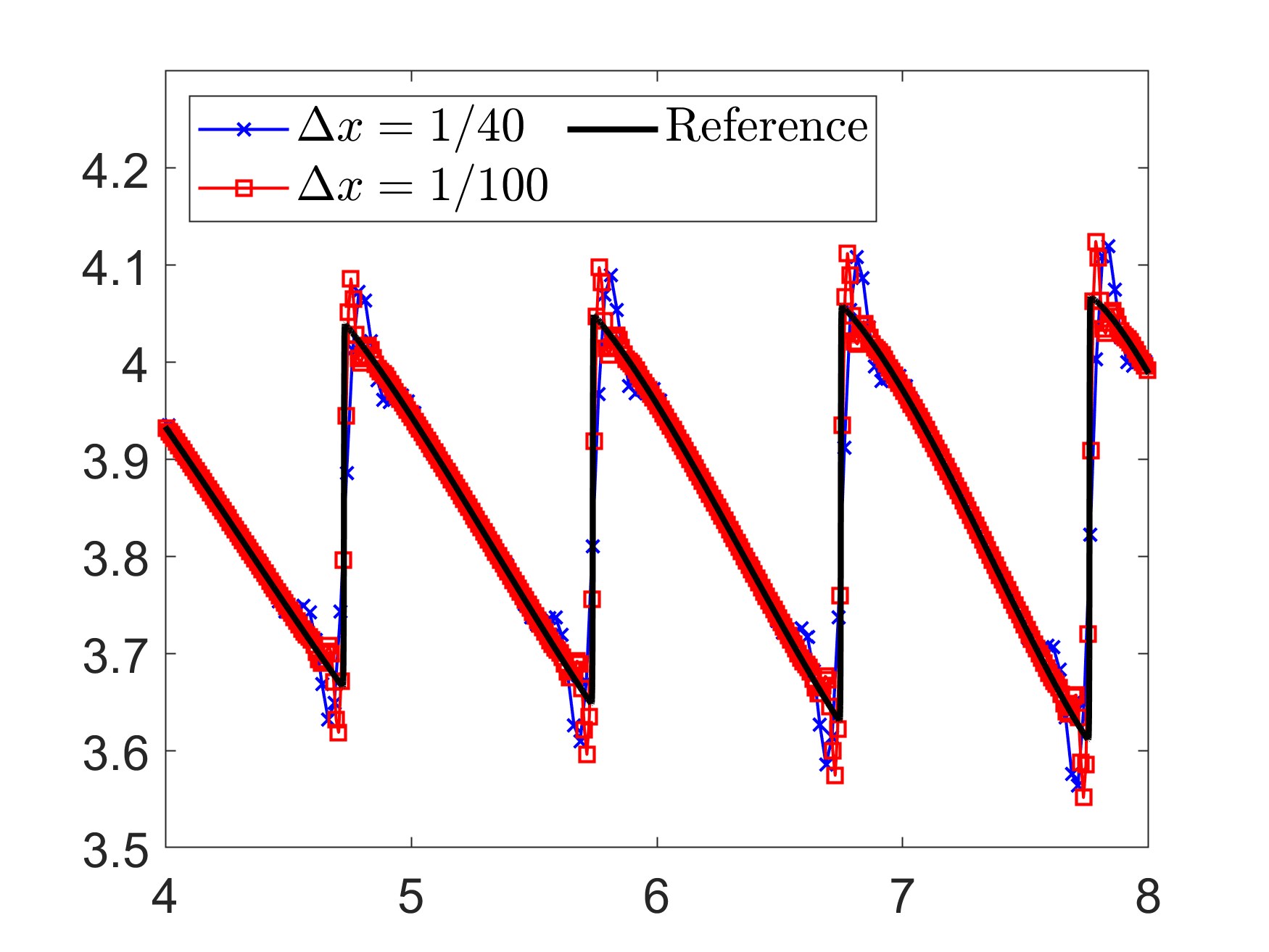}}
\caption{\sf Example 1: Density $\rho$ computed with $\dx=1/40$ and $\dx=1/100$ by the MODIFIED scheme (left) and zoom at
$x\in[4,8]$ (right).\label{fig1c}}
\end{figure}

We compute the numerical solutions on uniform meshes with $\dx=1/40$ and $1/100$ until the final time $t=5$ by the MODIFIED scheme with the
adaption constant $\texttt{C}_1=0.01$. We present the obtained densities in Figure \ref{fig1a}, where one can see that compared with the NEW
and OLD schemes (see Figures \ref{fig4} and \ref{fig4a}), the MODIFIED one produces ``stair-like'' structures, which are especially
pronounced in a finer mesh solution. The appearance of these nonphysical structures suggests that the value $\texttt{C}_1=0.01$ is too small
and hence the identified ``rough'' areas, where the overcompressive limiter is applied, are too wide. We therefore take a larger value
$\texttt{C}_1=0.1$ and repeat the same computations on the same two meshes. For this value of $\texttt{C}_1$, no ``stair-like'' structures
are formed even when the mesh is refined. However, the obtained densities, plotted in Figure \ref{fig1c}, contain oscillations in the
vicinities of shock waves. This suggests that now the identified smooth area is too wide and the quasi-linear fifth-order scheme is used
near the shocks. It is noticeable that the magnitude of oscillations near the shocks does not decay when the mesh is refined.

We emphasize that the need to employ a slightly more complicated, but far more accurate and robust, NEW adaptive scheme is justified by the
failure of a simpler MODIFIED adaptive scheme in this example. Therefore, we will only compare the NEW and OLD schemes' performance in the
rest of \S\ref{sec4}.

\paragraph{Example 2---Titarev-Toro Problem.} In the second example taken from \cite{Toro2005} (see also \cite{Shu88,Toro2005a}), we
consider the shock-entropy wave interaction problem with the following initial conditions:
\begin{equation*}
(\rho, u,p)\Big|_{(x,0)}=\begin{cases}(1.51695,0.523346,1.805),&x<-4.5,\\(1+0.1\sin(20x),0,1),&x>-4.5,\end{cases}
\end{equation*}
which correspond to a forward-facing shock wave of Mach number $1.1$ interacting with high-frequency density perturbations, that is, as the
shock wave moves, the perturbations spread ahead. We set the free boundary condition at both ends of the computational domain $[-5,5]$.

\begin{figure}[ht!]
\centerline{\includegraphics[trim=0.9cm 0.4cm 1.2cm 0.6cm, clip, width=6.0cm]{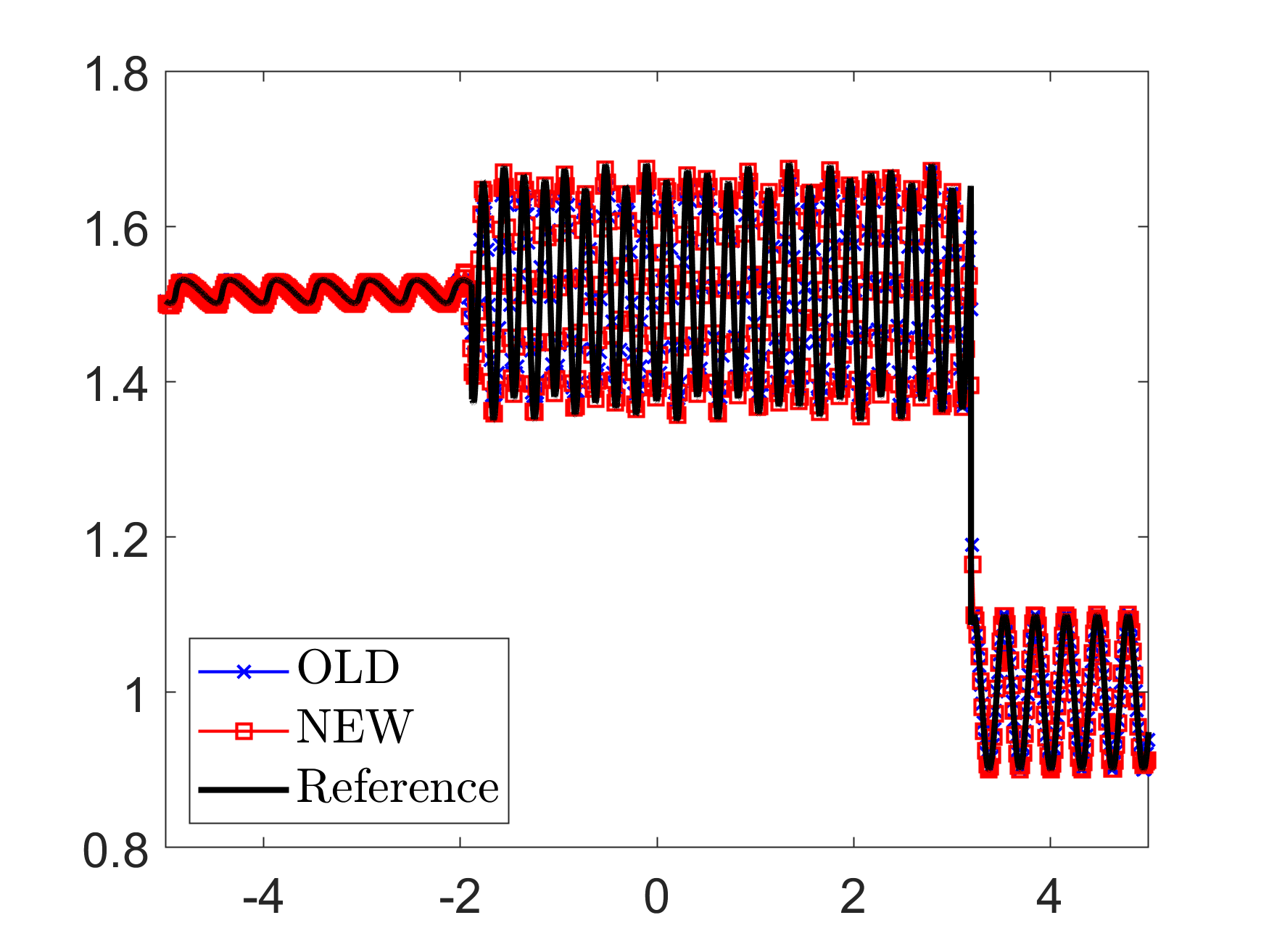}\hspace{1cm}
            \includegraphics[trim=0.9cm 0.4cm 1.2cm 0.6cm, clip, width=6.0cm]{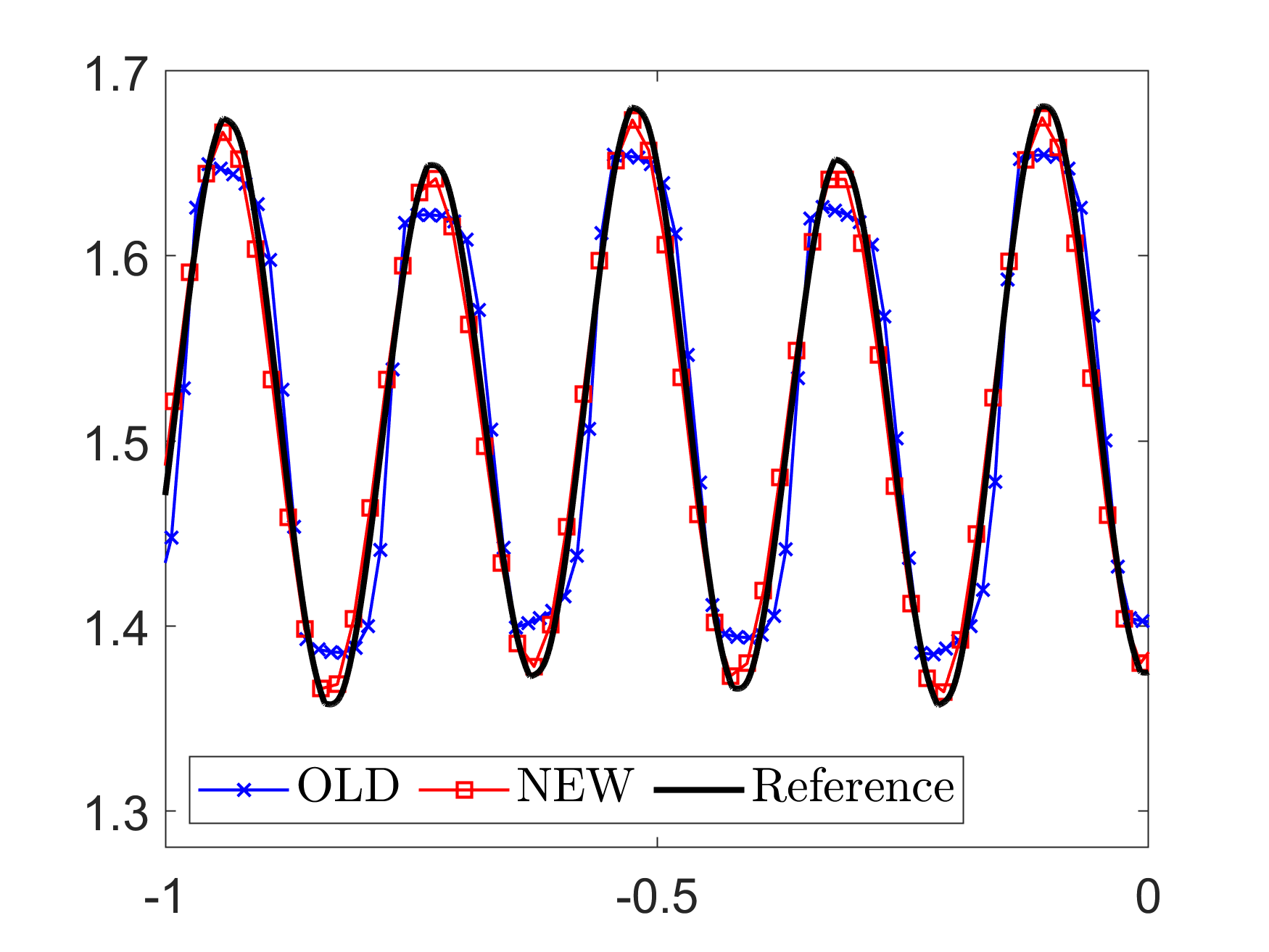}}
\caption{\sf Example 2: Density $\rho$ computed by the NEW and OLD  schemes (left) and zoom at $x\in[-1,0]$ (right).\label{fig3}}
\end{figure}

We compute the numerical solution until the final time $t=5$ by the NEW (with the adaption constants $\texttt{C}_1=0.02$ and
$\texttt{C}_2=0.3$) and OLD (with the adaption constant $\texttt{C}_1=0.01$) schemes on a uniform mesh with $\dx=1/60$ . The obtained
numerical results are presented in Figure \ref{fig3} along with the reference solution computed on a much finer mesh with $\dx=1/1600$,
where one can see that the NEW scheme is substantially more accurate than the OLD one.

\paragraph*{Example 3---Blast Wave Problem.} In the last 1-D example, we consider a strong shocks interaction problem from \cite{Woodward88}
with the initial data,
\begin{equation*}
(\rho,u,p)\Big|_{(x,0)}=\begin{cases}(1,0,1000),&x<0.1,\\(1,0,0.01),&0.1\le x\le0.9,\\(1,0,100),&x>0.9,\end{cases}
\end{equation*}
prescribed in the computational domain $[0,1]$ subject to the solid wall boundary conditions.

\begin{figure}[ht!]
\centerline{\includegraphics[trim=1.3cm 0.4cm 0.8cm 0.6cm, clip, width=6.0cm]{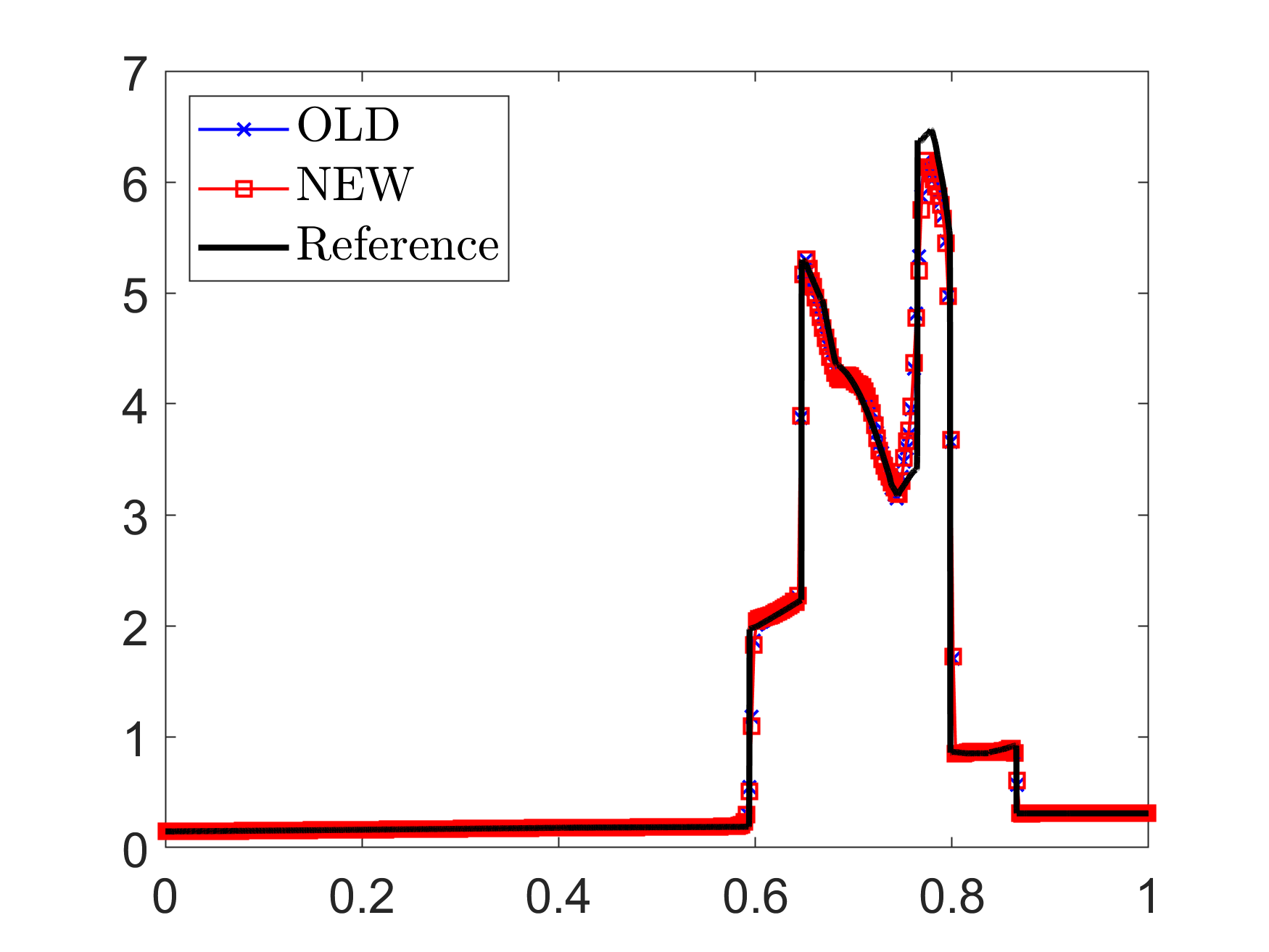}\hspace*{1.0cm}
            \includegraphics[trim=1.3cm 0.4cm 0.8cm 0.6cm, clip, width=6.0cm]{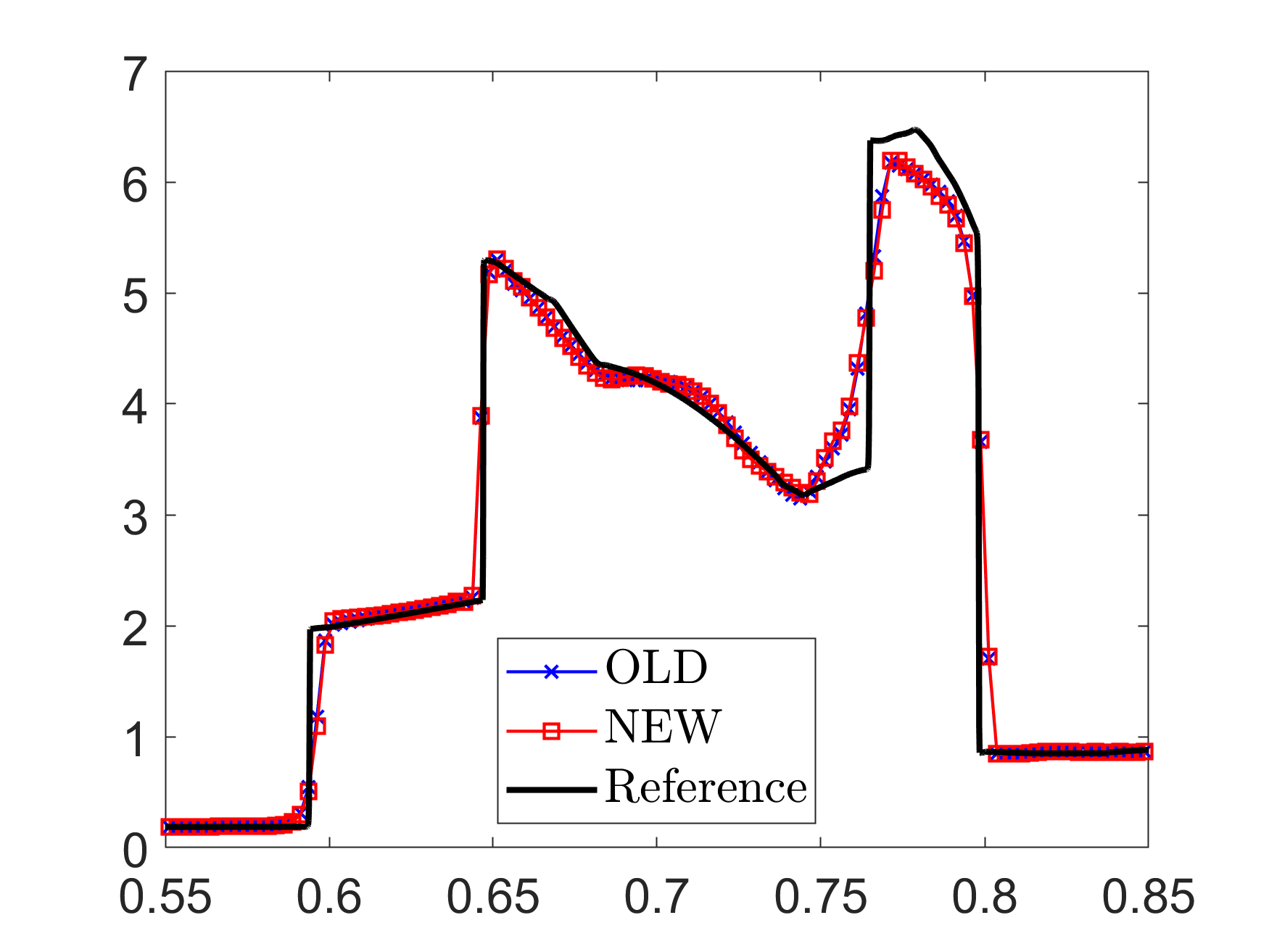}}
\caption{\sf Example 3: Density $\rho$ computed by the NEW and OLD schemes (left) and zoom at $x\in[0.55,0.85]$.\label{fig5}}
\end{figure}

We compute the numerical solutions until the final time $t=0.038$ by the NEW (with the adaption constants $\texttt{C}_1=0.02$ and
$\texttt{C}_2=0.3$) and OLD (with the adaption constant $\texttt{C}_1=0.01$) schemes on a uniform mesh with $\dx=1/400$ and present the
obtained numerical results in Figure \ref{fig5} together with the reference solution computed on a much finer mesh with $\dx=1/4000$. It can
be observed that both the NEW and OLD schemes are capable of achieving a superb resolution of the contact wave located at about $x=0.6$.

\subsection{2-D Examples}
We now turn to the 2-D numerical examples. In Examples 4--6, we will consider three different 2-D Riemann problems from \cite{Kurganov02};
see also \cite{Schulz93,Schulz93a,Zheng01}.

In all of the examples below, we will show the computed densities and the corresponding ``rough'' areas detected by the SIs. In the plots
obtained by the NEW scheme, we also highlight in red color the detected ``contact'' subareas of the ``rough'' areas.

\paragraph*{Example 4---2-D Riemann Problem (Configuration 3).} In this example, the initial conditions,
\begin{equation*}
(\rho,u,v,p)\Big|_{(x,y,0)}=\begin{cases}(1.5,0,0,1.5),&x>1,~y>1,\\(0.5323,1.206,0,0.3),&x<1,~y>1,\\(0.138,1.206,1.206,0.029),&x<1,~y<1,\\
(0.5323,0,1.206,0.3),&x>1,~y<1,\end{cases}
\end{equation*}
are prescribed in the computational domain $[0,1.2]\times[0,1.2]$ subject to the free boundary conditions.

\begin{figure}[ht!]
\centerline{\includegraphics[trim=3.6cm 3.3cm 1.2cm 2.5cm, clip, width=12.6cm]{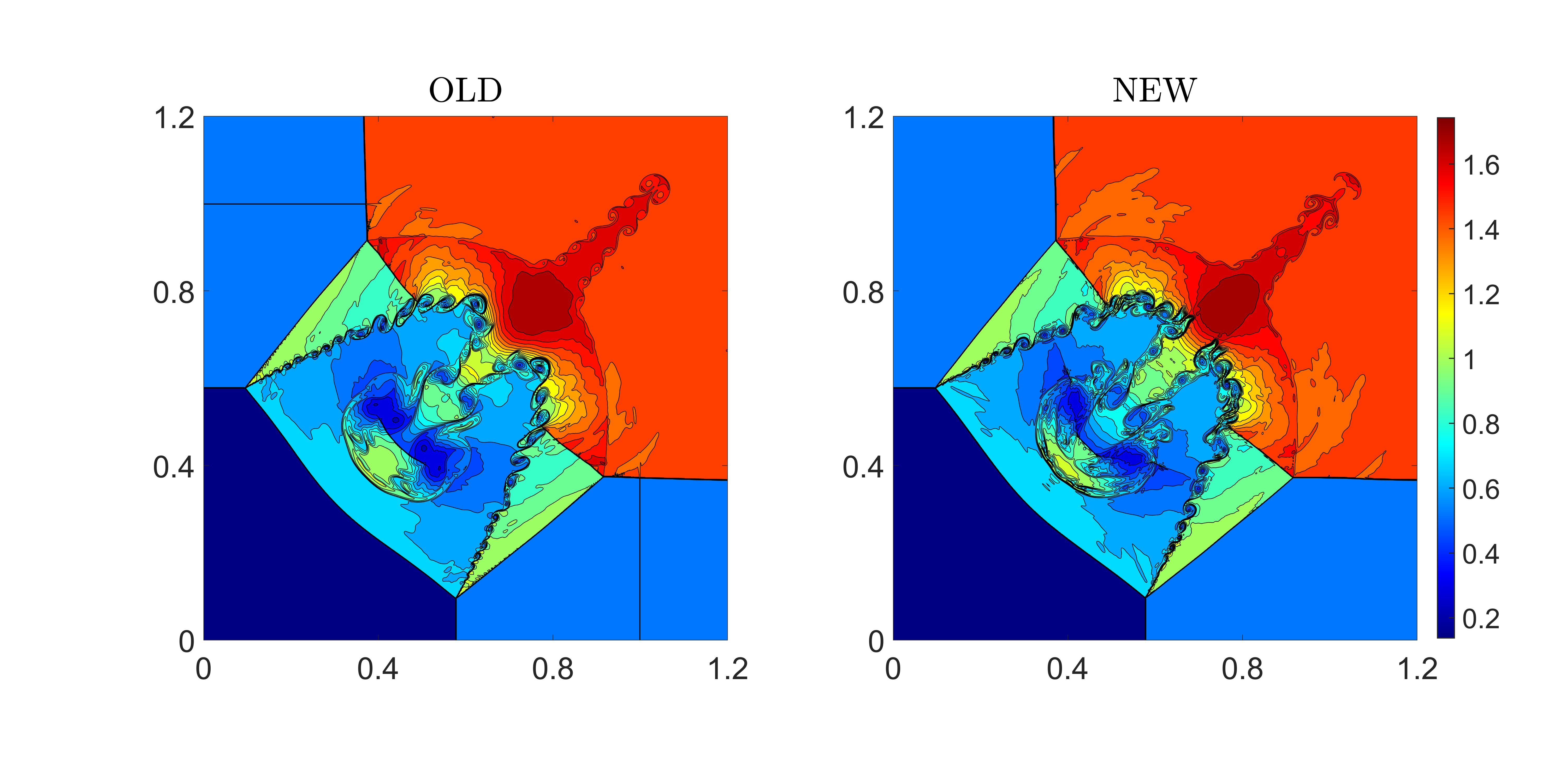}}
\caption{\sf Example 4: Density $\rho$ computed by the OLD (left) and NEW (right) schemes.\label{fig6a}}
\end{figure}

\begin{figure}[ht!]
\centerline{\hspace*{-0.5cm}\includegraphics[trim=2.1cm 0.3cm 2.2cm 0.1cm, clip, width=5.2cm]{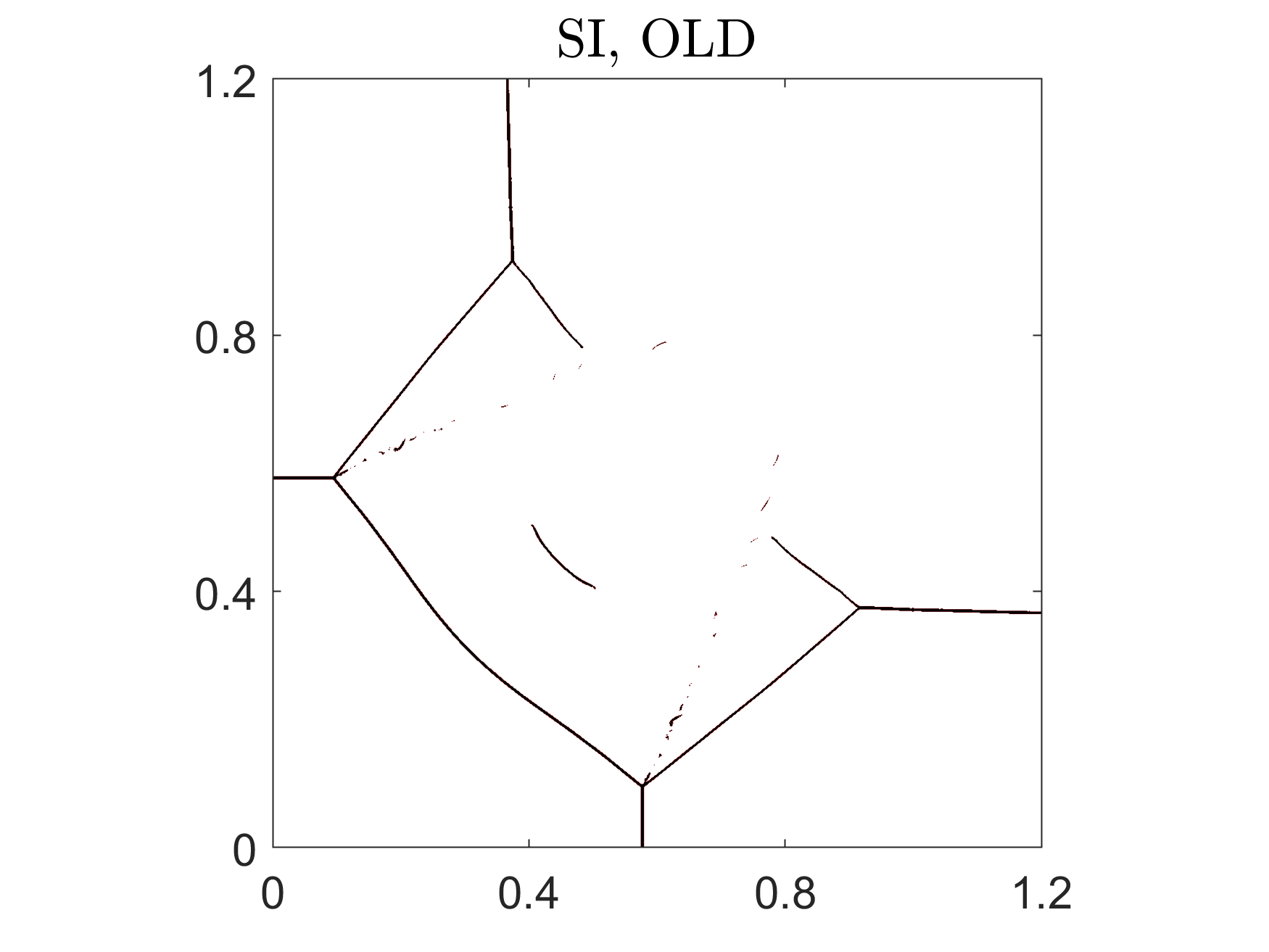}\hspace*{0.8cm}
                            \includegraphics[trim=2.1cm 0.3cm 2.2cm 0.1cm, clip, width=5.2cm]{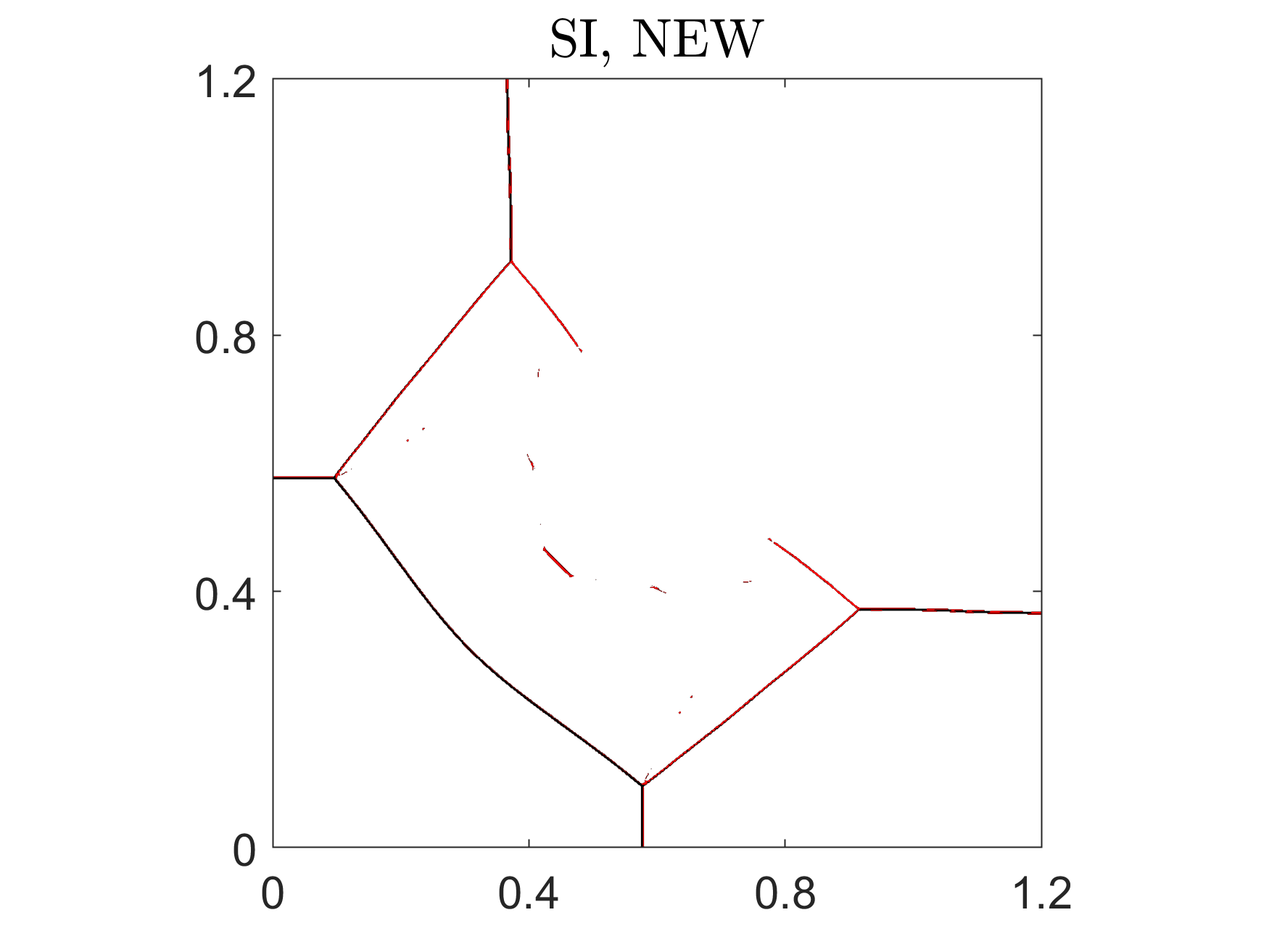}}
\caption{\sf Example 4: The ``rough'' areas detected by the OLD (left) and NEW (right) schemes at the final time-step.\label{fig6b}}
\end{figure}

We compute the numerical solution until the final time $t=1$ by the NEW (with the adaption constants $\texttt{C}_1=0.08$ and
$\texttt{C}_2=0.2$) and OLD (with the adaption constant $\texttt{C}_1=0.08$) schemes on a uniform mesh with $\dx=\dy=3/2500$. The obtained
results are plotted in Figure \ref{fig6a}, where one can see that the NEW scheme slightly outperforms the OLD one in capturing a sideband
instability of the jet in the zones of strong along-jet velocity shear and the instability along the jets neck. In Figure \ref{fig6b}, we
show the detected ``rough'' areas (for both the NEW and OLD schemes) and their ``contact'' subareas (for the NEW scheme). As one can see,
when the NEW scheme is used, the SBM and Minmod2 limiters are implemented only in a small part of the computational domain.

\paragraph*{Example 5---2-D Riemann Problem (Configuration 6).} In this example, the initial conditions,
\begin{equation*}
(\rho,u,v,p)\Big|_{(x,y,0)}=\begin{cases}(1,0.75,-0.5,1),&x>0.5,~y>0.5,\\(2,0.75,0.5,1),&x<0.5,~y>0.5,\\(1,-0.75,0.5,1),&x<0.5,~y<0.5,\\
(3,-0.75,-0.5,1),&x>0.5,~y<0.5,\end{cases}
\end{equation*}
are prescribed in the computational domain $[0,1]\times[0,1]$ subject to the free boundary conditions.

\begin{figure}[ht!]
\centerline{\includegraphics[trim=5.3cm 6.7cm 2.4cm 5.7cm, clip, width=14.5cm]{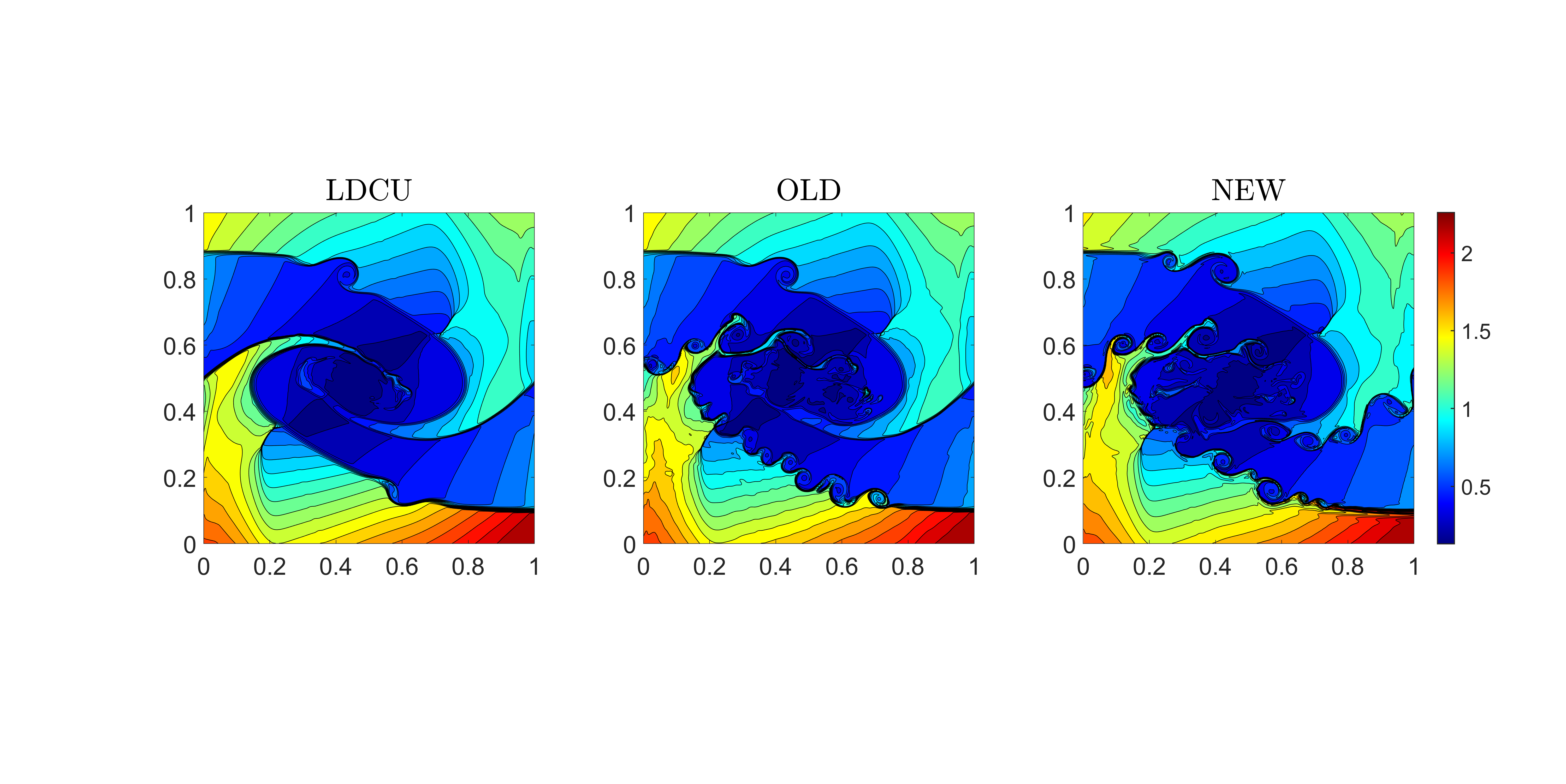}}
\caption{\sf Example 5: Density $\rho$ computed by the LDCU (left), OLD (middle), and NEW (right) schemes.\label{fig16a}}
\end{figure}

\begin{figure}[ht!]
\centerline{\hspace*{4.1cm}\includegraphics[trim=2.1cm 0.3cm 2.2cm 0.1cm, clip, width=4.1cm]{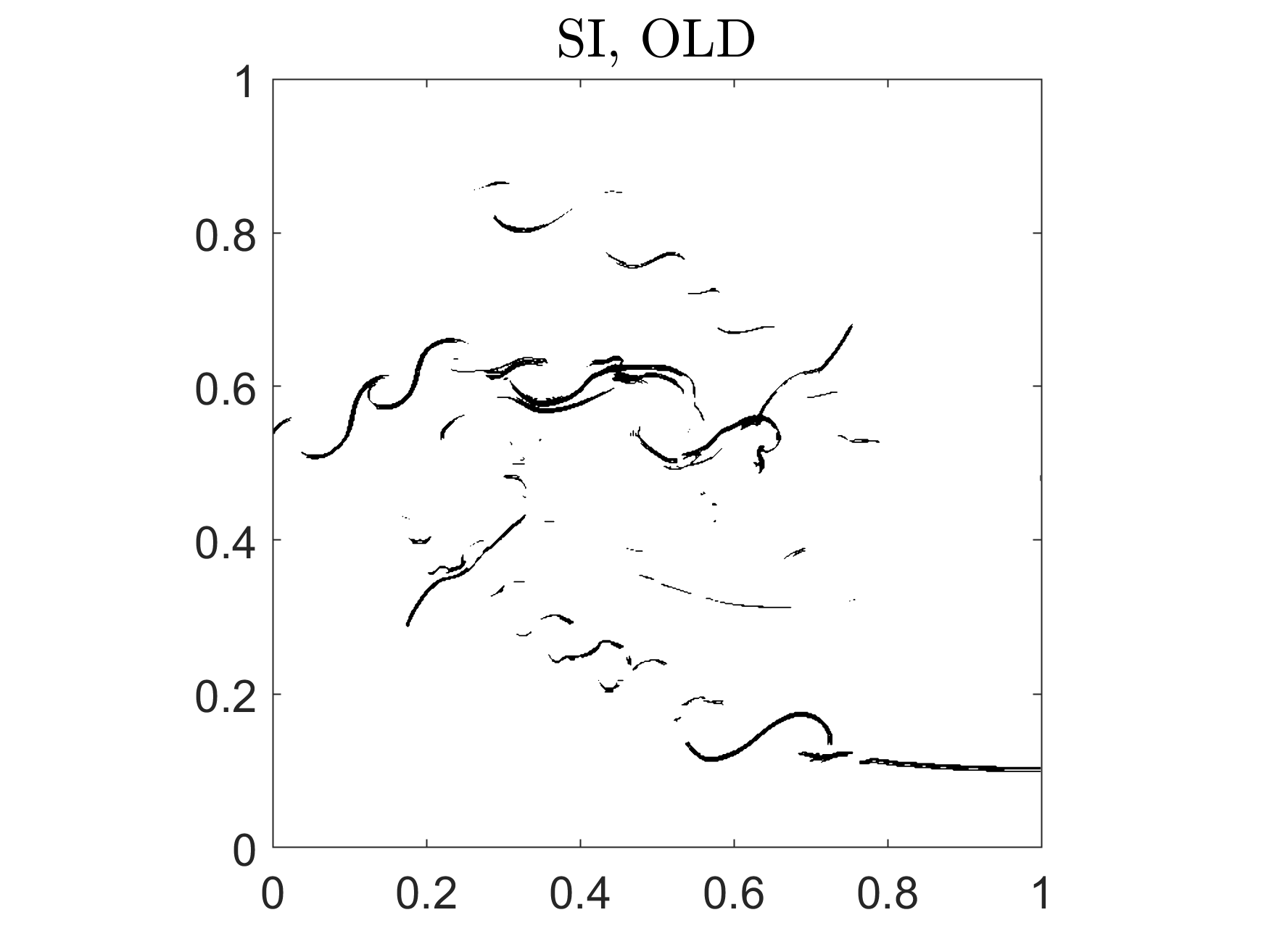}\hspace*{0.7cm}
                           \includegraphics[trim=2.1cm 0.3cm 2.2cm 0.1cm, clip, width=4.1cm]{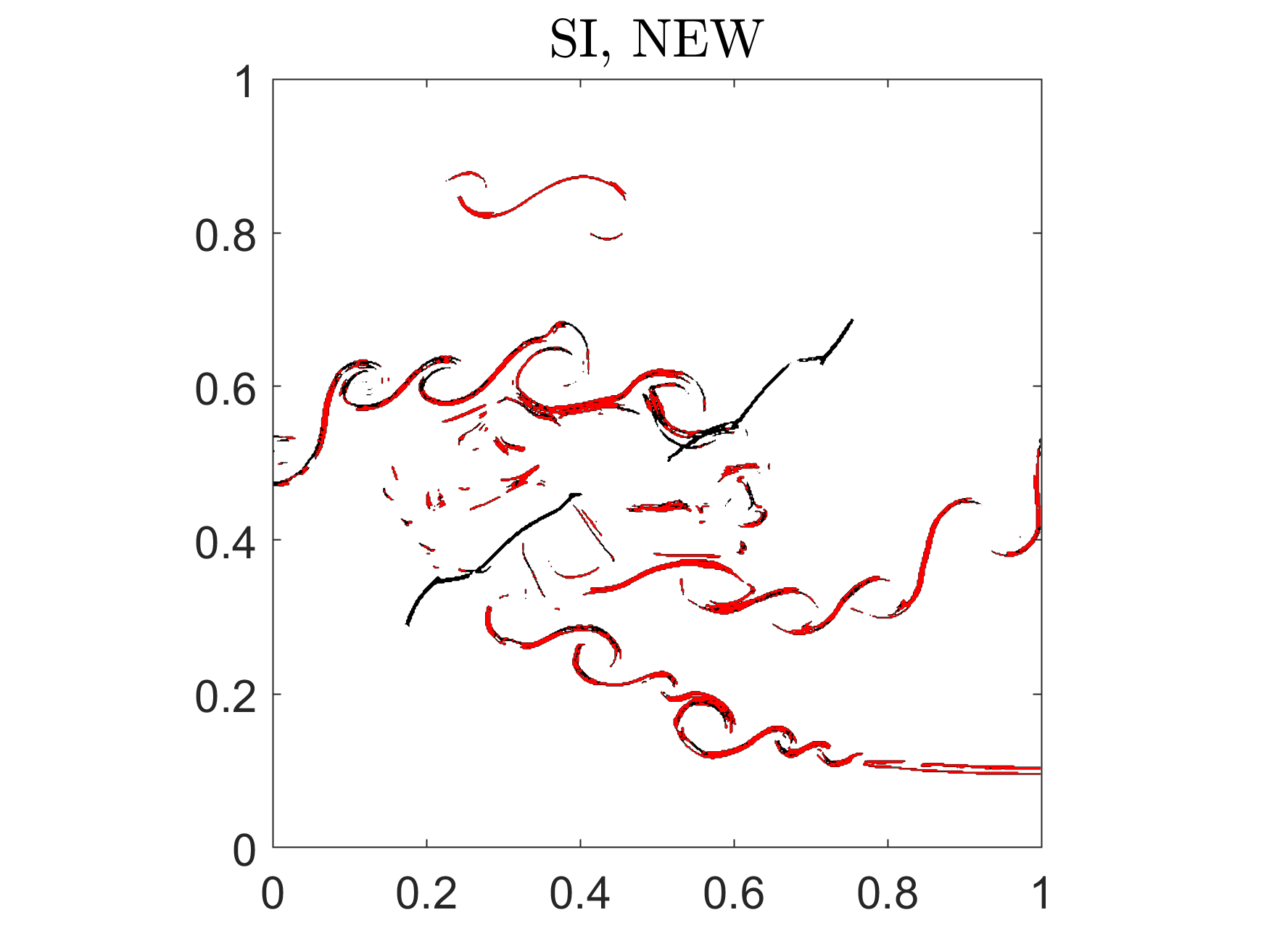}}
\caption{\sf Example 5: The ``rough'' areas detected by the OLD (left) and NEW (right) schemes at the final time-step.\label{fig16b}}
\end{figure}

We compute the numerical solution until the final time $t=1$ by the second-order LDCU (implemented with the Minmod2 limiter), NEW (with the
adaption constants $\texttt{C}_1=0.06$ and $\texttt{C}_2=0.004$), and OLD (with the adaption constant $\texttt{C}_1=0.1$) schemes on a
uniform mesh with $\dx=\dy=1/600$, and plot the obtained results in Figure \ref{fig16a}, where one can see that the adaptive schemes
clearly outperform the LDCU one. At the same time, the NEW scheme is capable of capturing much more complicated vortex structures compared
with the OLD scheme, which demonstrates higher resolution and lower numerical dissipation of the NEW scheme. Note that the ``rough'' areas
identified by both the NEW and OLD schemes are quite similar; see Figure \ref{fig16b}.

\paragraph*{Example 6---2-D Riemann Problem (Configuration 12).} In this example, the initial conditions,
\begin{equation*}
(\rho,u,v,p)\Big|_{(x,y,0)}=\begin{cases}(0.5313,0,0,0.4),&x>0.5,~y>0.5,\\(1,0.7276,0,1),&x<0.5,~y>0.5,\\(0.8,0,0,1),&x<0.5,~y<0.5,\\
(1,0,0.7276,1),&x>0.5,~y<0.5,\end{cases}
\end{equation*}
are prescribed in the computational domain $[0,0.6]\times[0,0.6]$ subject to the free boundary conditions.

\begin{figure}[ht!]
\centerline{\includegraphics[trim=5.3cm 6.7cm 2.4cm 5.7cm, clip, width=14.5cm]{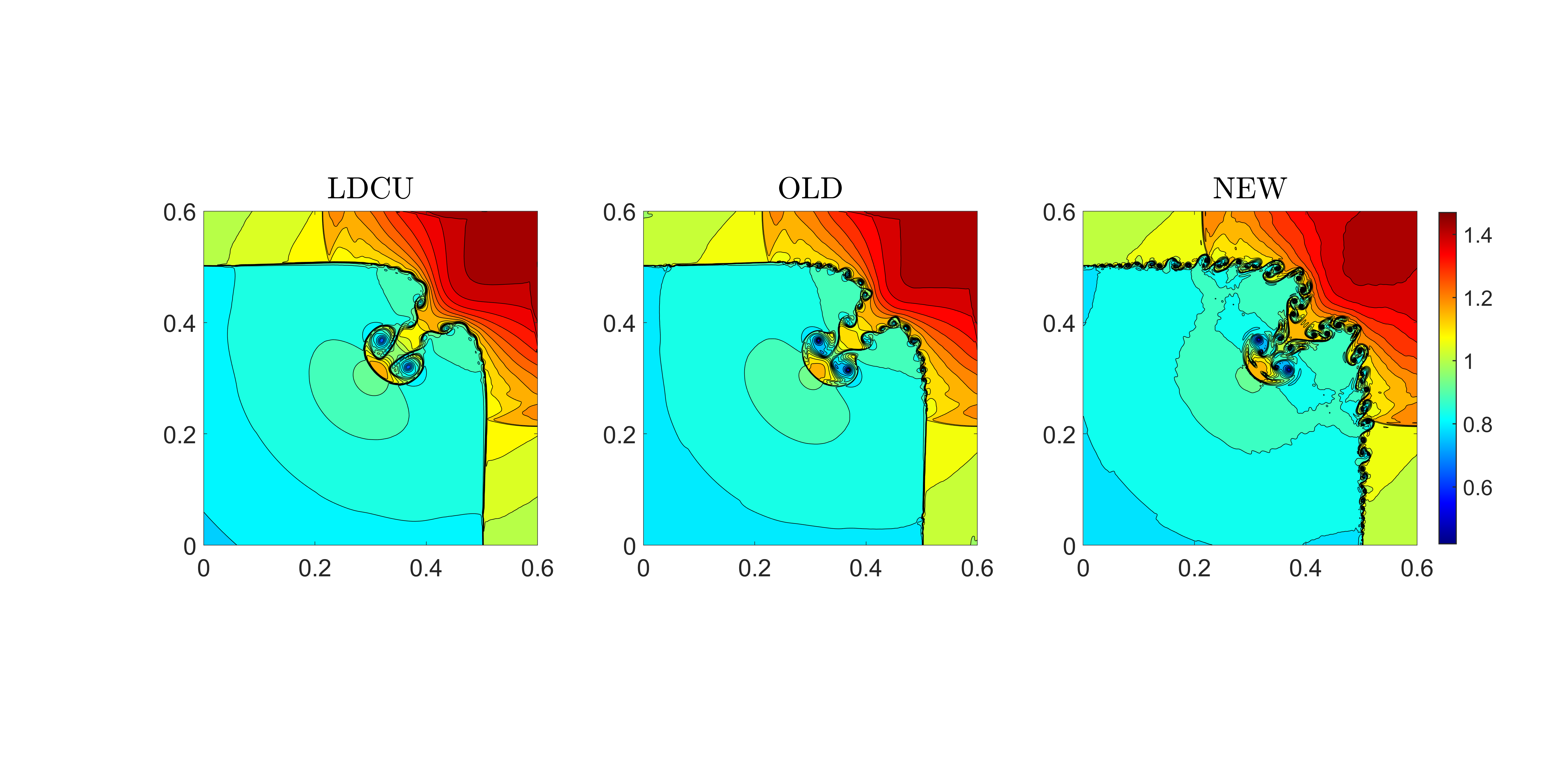}}
\caption{\sf Example 6: Density $\rho$ computed by the LDCU (left), OLD (middle), and NEW (right) schemes.\label{fig17a}}
\end{figure}

\begin{figure}[ht!]
\centerline{\hspace*{4.1cm}\includegraphics[trim=2.1cm 0.3cm 2.2cm 0.1cm, clip, width=4.1cm]{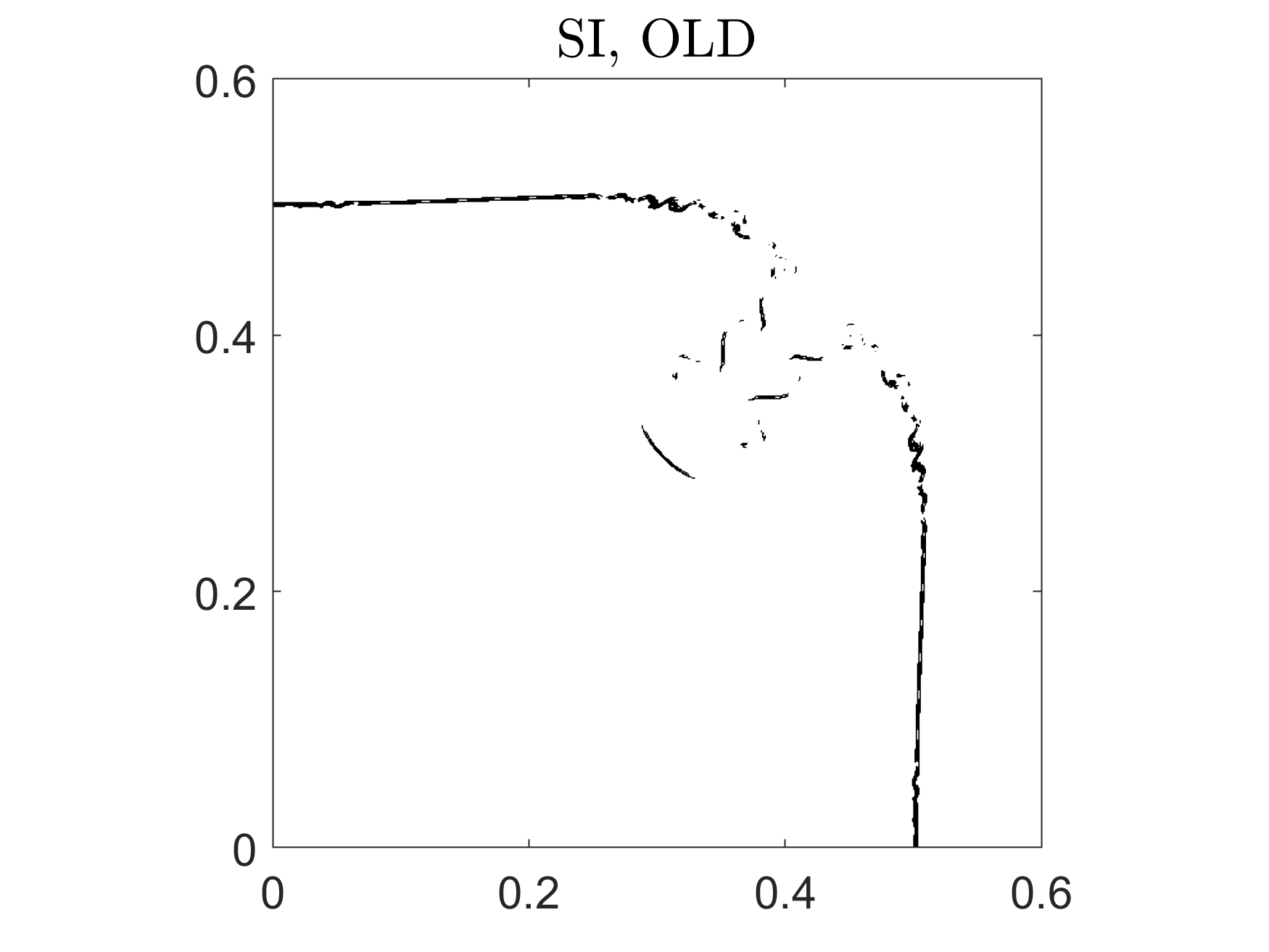}\hspace*{0.7cm}
                           \includegraphics[trim=2.1cm 0.3cm 2.2cm 0.1cm, clip, width=4.1cm]{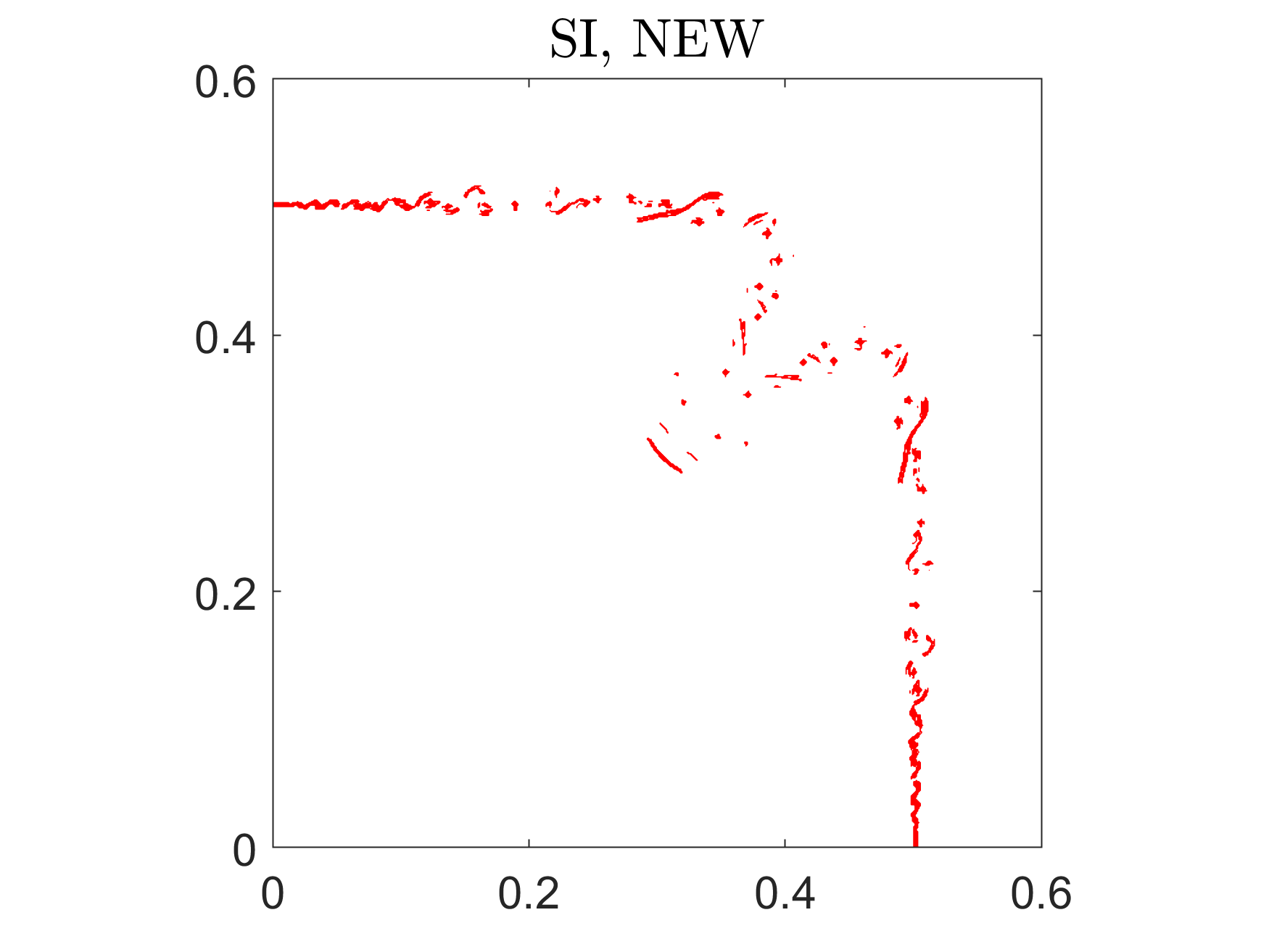}}
\caption{\sf Example 6: The ``rough'' areas detected by the OLD (left) and NEW (right) schemes at the final time-step.\label{fig17b}}
\end{figure}

We compute the numerical solution until the final time $t=0.5$ by the second-order LDCU (implemented with the Minmod2 limiter), NEW (with
the adaption constants $\texttt{C}_1=0.03$ and $\texttt{C}_2=0.2$), and OLD (with the adaption constant $\texttt{C}_1=0.03$) schemes on a
uniform mesh with $\dx=\dy=1/1000$, and plot the obtained results in Figure \ref{fig17a}. As one can see, the NEW scheme resolves
substantially more vortices, arising along the unstable contact surfaces, than both the OLD and LDCU schemes. As in the previous examples,
the ``rough'' areas detected by the NEW and OLD schemes are quite similar; see Figure \ref{fig17b}.

\paragraph{Example 7---Implosion Problem.} In this example, we consider the implosion problem taken from \cite{Liska03}. The initial
conditions,
\begin{equation*}
(\rho,u,v,p)\Big|_{(x,y,0)}=\begin{cases}(0.125,0,0,0.14),&|x|+|y|<0.15,\\(1,0,0,1),&\mbox{otherwise},\end{cases}
\end{equation*}
are prescribed in the computational domain $[0,0.5]\times[0,0.5]$ with the solid wall boundary conditions. In this setting, a jet forms near
the origin and propagates along the diagonal $y=x$ direction, and the position of the jet indicates the amount of numerical dissipation
present in studied numerical schemes. In general, schemes containing large numerical dissipation may not resolve the jet at all or the jet
propagation velocity may be affected by the numerical diffusion.

\begin{figure}[ht!]
\centerline{\includegraphics[trim=3.6cm 3.3cm 1.2cm 2.5cm, clip, width=12.6cm]{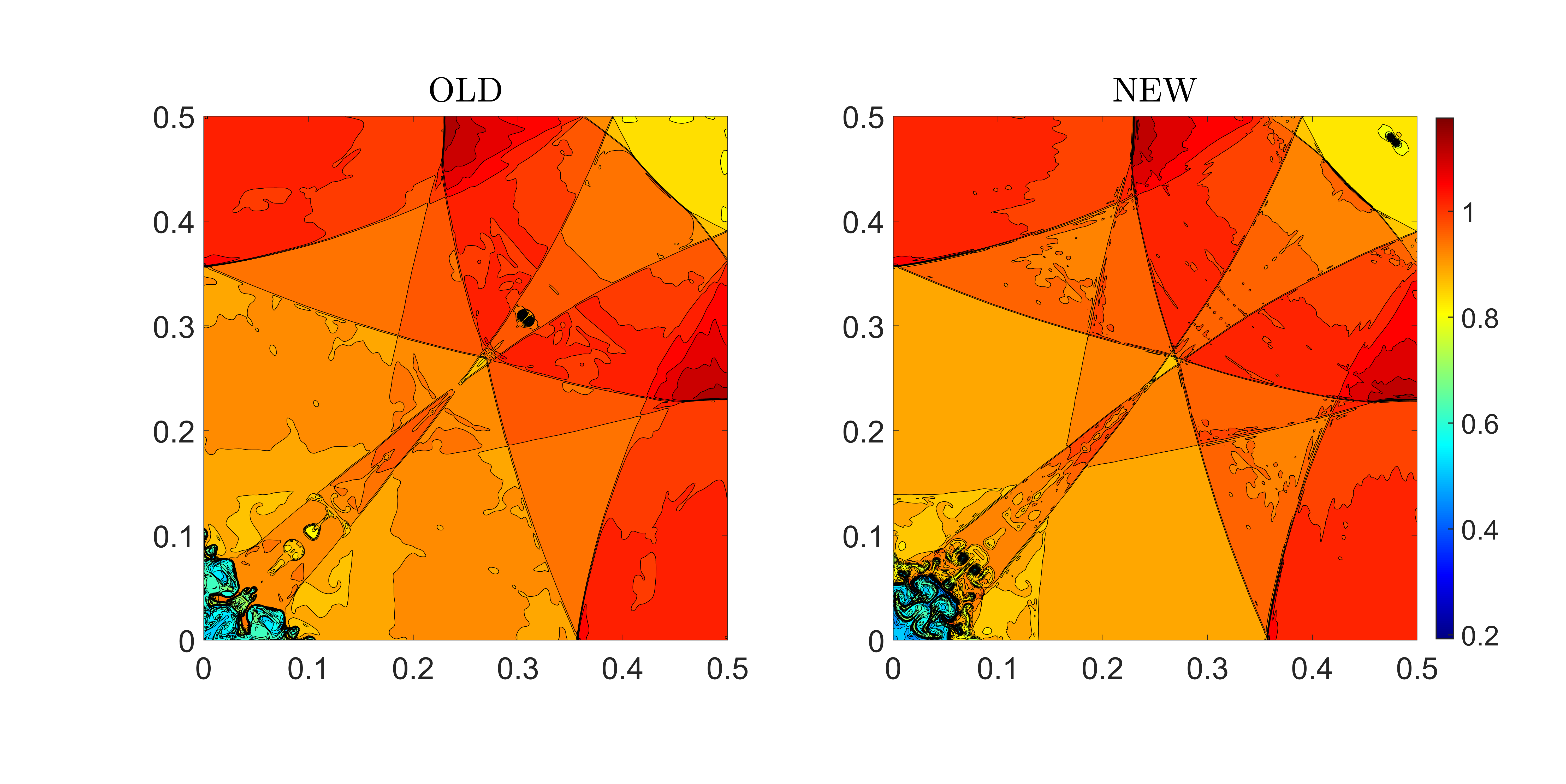}}
\caption{\sf Example 7: Density $\rho$ computed by the OLD (left) and NEW (right) schemes.\label{fig8a}}
\end{figure}

\begin{figure}[ht!]
\centerline{\hspace*{-0.5cm}\includegraphics[trim=2.1cm 0.3cm 2.2cm 0.1cm, clip, width=5.2cm]{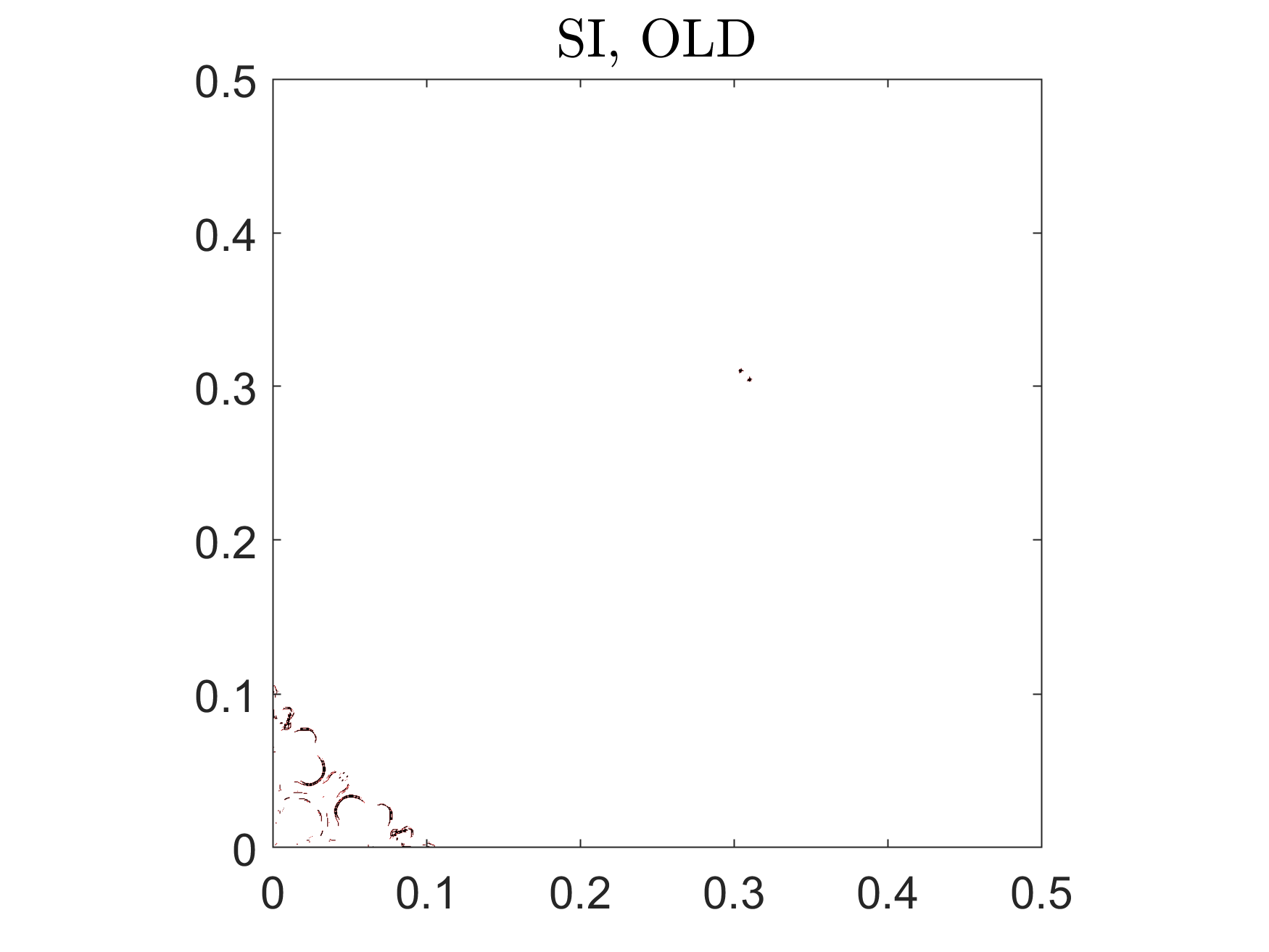}\hspace*{0.8cm}
                            \includegraphics[trim=2.1cm 0.3cm 2.2cm 0.1cm, clip, width=5.2cm]{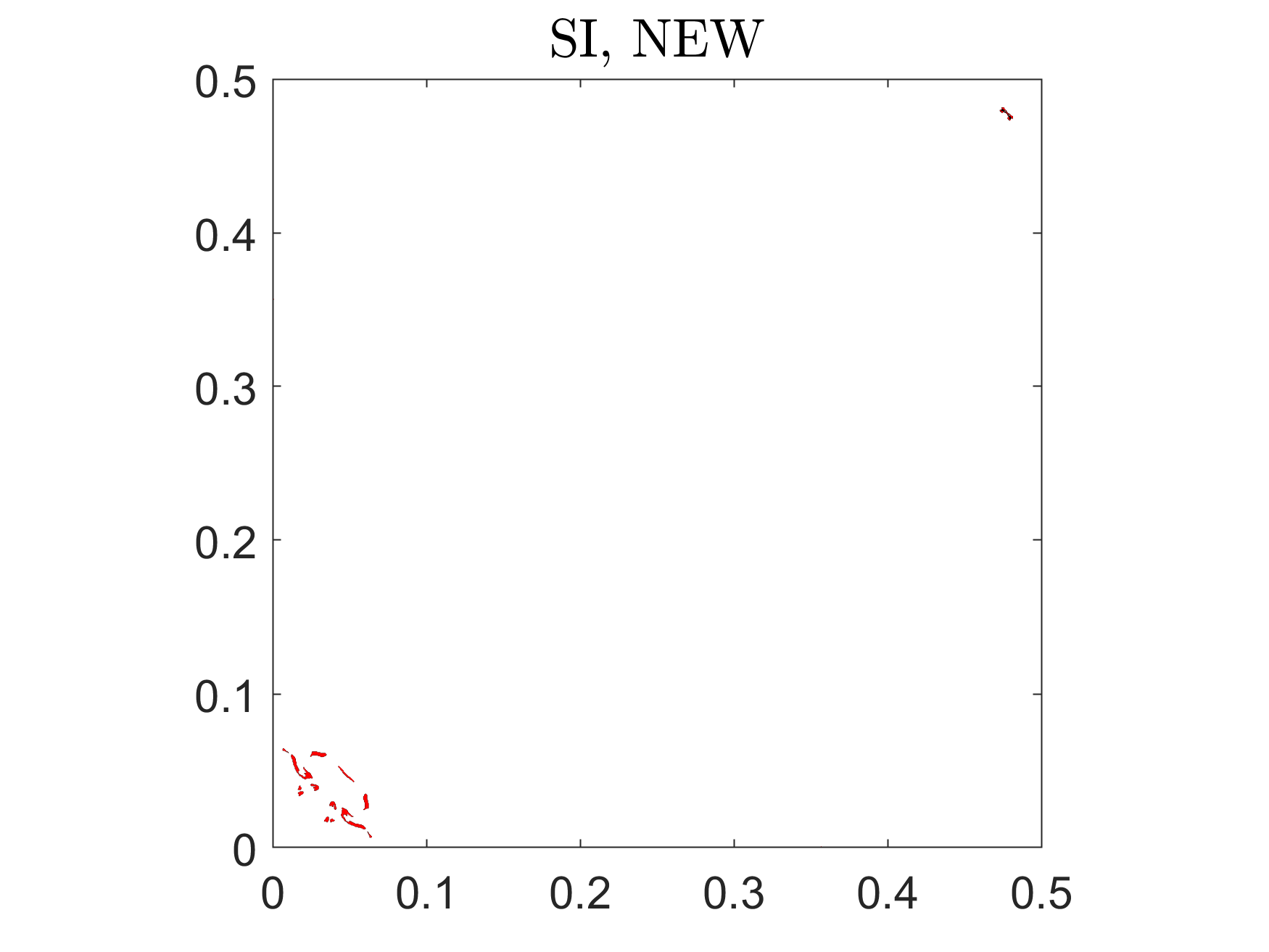}}
\caption{\sf Example 7: The ``rough'' areas detected by the OLD (left) and NEW (right) schemes at the final time-step.\label{fig8b}}
\end{figure}

We compute the numerical solution until the final time $t=2.5$ by the NEW (with the adaption constants $\texttt{C}_1=0.05$ and
$\texttt{C}_2=0.1$) and OLD (with the adaption constant $\texttt{C}_1=0.05$) schemes on a uniform mesh with $\dx=\dy=1/2000$ and present the
obtained numerical results in Figure \ref{fig8a}. As one can observe, compared with the OLD scheme, the jet propagates much further in the
diagonal direction in the numerical result computed by the NEW scheme, which demonstrates that the NEW scheme is less dissipative. In Figure
\ref{fig8b}, we show ``rough'' areas detected by the NEW and OLD schemes. One can see that when the NEW scheme is used, a sharper SBM
limiter is implemented only in a very small part of the computational domain, which is substantially smaller than the one detected as
``rough'' by the OLD scheme.

Next, we measure the CPU times consumed by the NEW and OLD schemes. The obtained results show that the CPU time consumed by the OLD scheme
is about $122\%$ larger than the CPU time consumed by the NEW scheme. This demonstrates that like in the 1-D case, the 2-D NEW scheme is
substantially more efficient than the OLD one though the efficiency gain in the 2-D case is not as big as in the 1-D one (see Example 1).

\paragraph{Example 8---Double Mach Reflection Problem.} In this example, we consider the double Mach reflection problem taken from
\cite{Woodward88}. The initial conditions,
\begin{equation*}
(\rho,u,v,p)\Big|_{(x,y,0)}=\begin{cases}(\rho_{\rm in},u_{\rm in},v_{\rm in},p_{\rm in}),&x<\dfrac{1}{6}+\dfrac{y}{\sqrt{3}},\\
(\rho_{\rm out},u_{\rm out},v_{\rm out},p_{\rm out}),&\mbox{otherwise},\end{cases}
\end{equation*}
with the inflow and outflow values,
\begin{equation*}
(\rho_{\rm in},u_{\rm in},v_{\rm in},p_{\rm in}):=\Big(8,\dfrac{17\sqrt{3}}{4},-\dfrac{17}{4},116.5\Big)\quad\mbox{and}\quad
(\rho_{\rm out},u_{\rm out},v_{\rm out},p_{\rm out}):=(1.4,0,0,1),
\end{equation*}
are prescribed in the computational domain $[0,4]\times[0,1]$. The supersonic Dirichlet boundary conditions are specified at $x=0$ and the
left part of the lower boundary $y=0$:
$$
(\rho,u,v,p)\Big|_{(0,y,t)}=(\rho_{\rm in},u_{\rm in},v_{\rm in},p_{\rm in}),~y\in[0,1],\qquad
(\rho,u,v,p)\Big|_{(x,0,t)}=(\rho_{\rm in},u_{\rm in},v_{\rm in},p_{\rm in}),~x\in\Big[0,\frac{1}{6}\Big],
$$
the free boundary conditions are imposed at $x=4$, and the solid wall boundary conditions are applied at the right part of the lower
boundary $y=0$ for $x\in\big[\frac{1}{6},4\big]$. Finally, at the upper boundary $y=1$, the exact Mach $10$ moving oblique shock is imposed,
namely, we set
\begin{equation*}
(\rho,u,v,p)\Big|_{(x,1,t)}=\begin{cases}(\rho_{\rm in},u_{\rm in},v_{\rm in},p_{\rm in}),&x<\dfrac{1}{6}+\dfrac{1+20t}{\sqrt{3}},\\
(\rho_{\rm out},u_{\rm out},v_{\rm out},p_{\rm out}),&\mbox{otherwise}.\end{cases}
\end{equation*}

This example is designed to assess the ability of the studied schemes to capture the complex shock and vortex structures generated when a
planar shock impinges on a wedge. In particular, it challenges the schemes to accurately resolve the progressive sharpening of the shock
that connects the contact surface with the transverse wave.

\begin{figure}[ht!]
\centerline{\includegraphics[trim=0.7cm 3.2cm 1.2cm 2.6cm, clip, width=7.3cm]{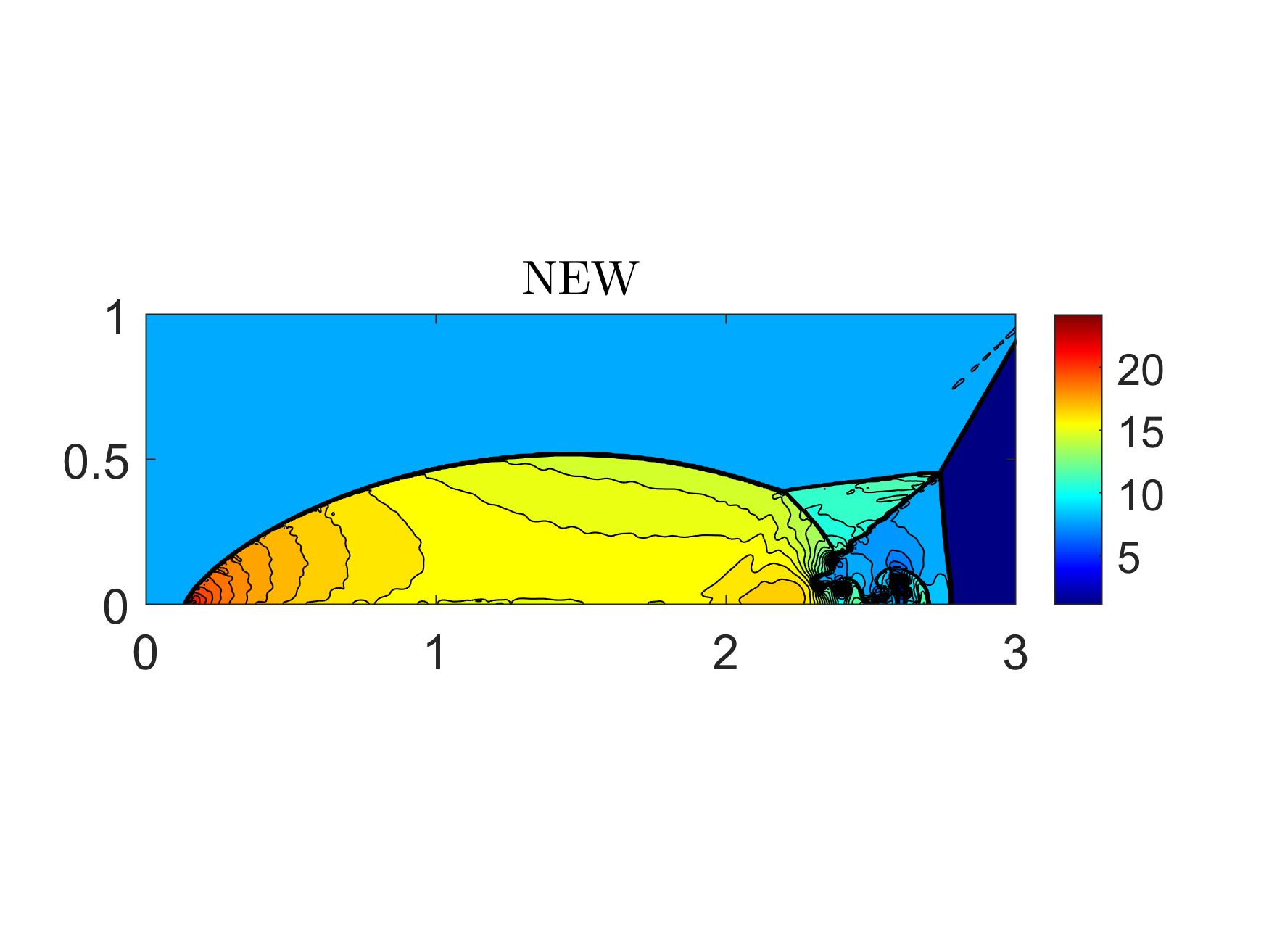}\hspace{0.5cm}
            \includegraphics[trim=0.7cm 2.9cm 1.2cm 2.6cm, clip, width=6.3cm]{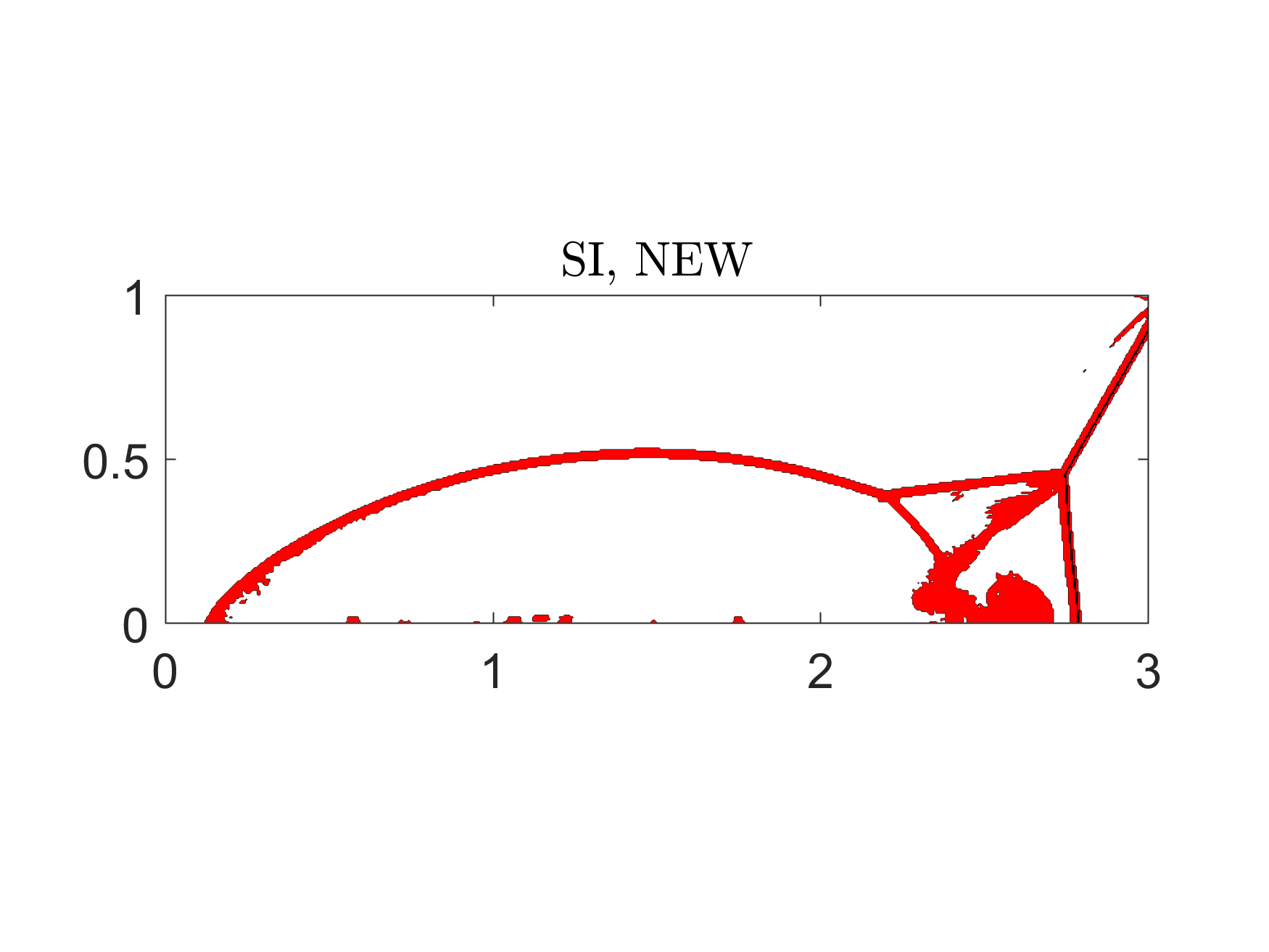}}
\centerline{\includegraphics[trim=0.7cm 3.2cm 1.2cm 2.6cm, clip, width=7.3cm]{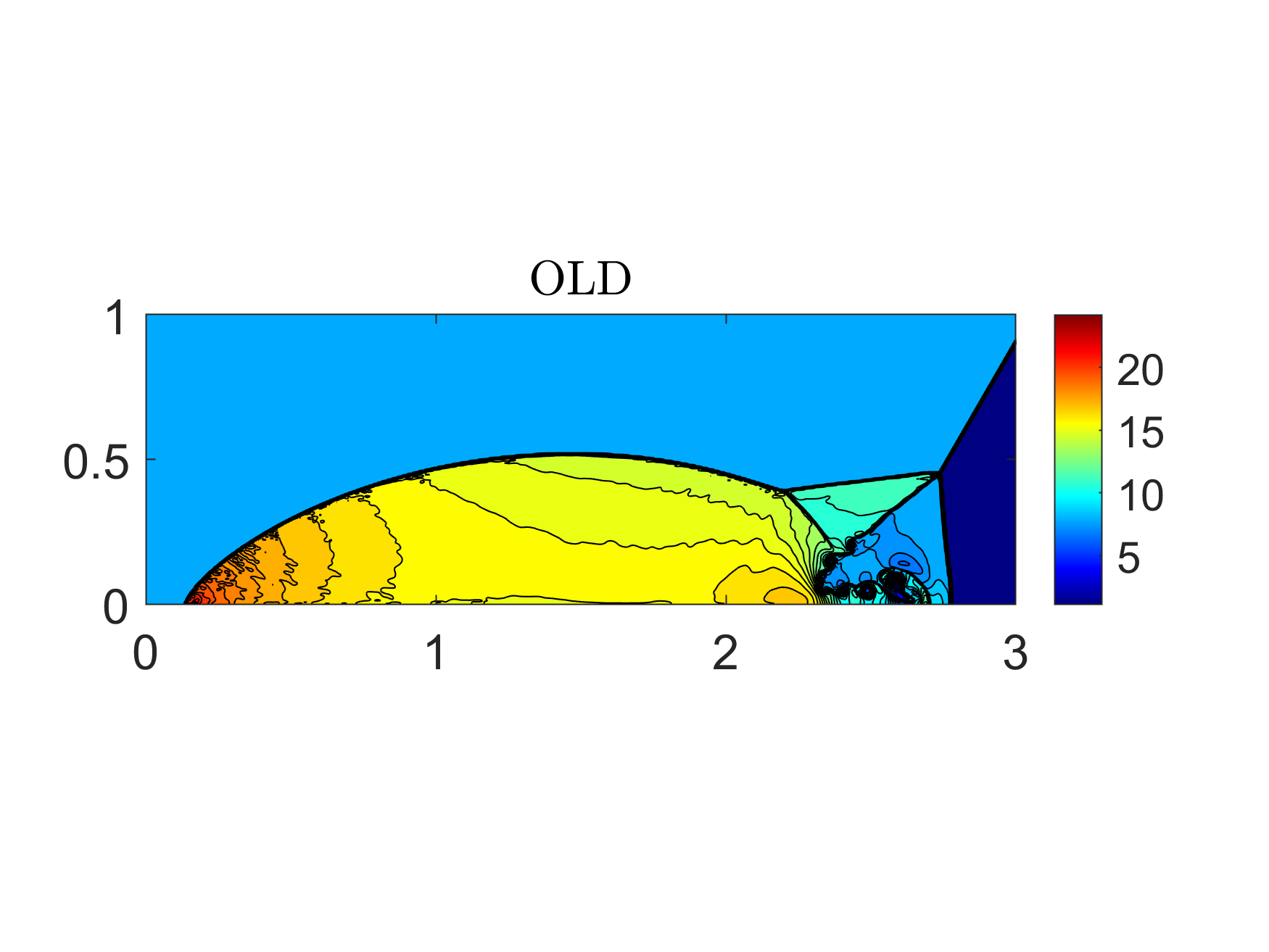}\hspace{0.5cm}
            \includegraphics[trim=0.7cm 2.9cm 1.2cm 2.6cm, clip, width=6.3cm]{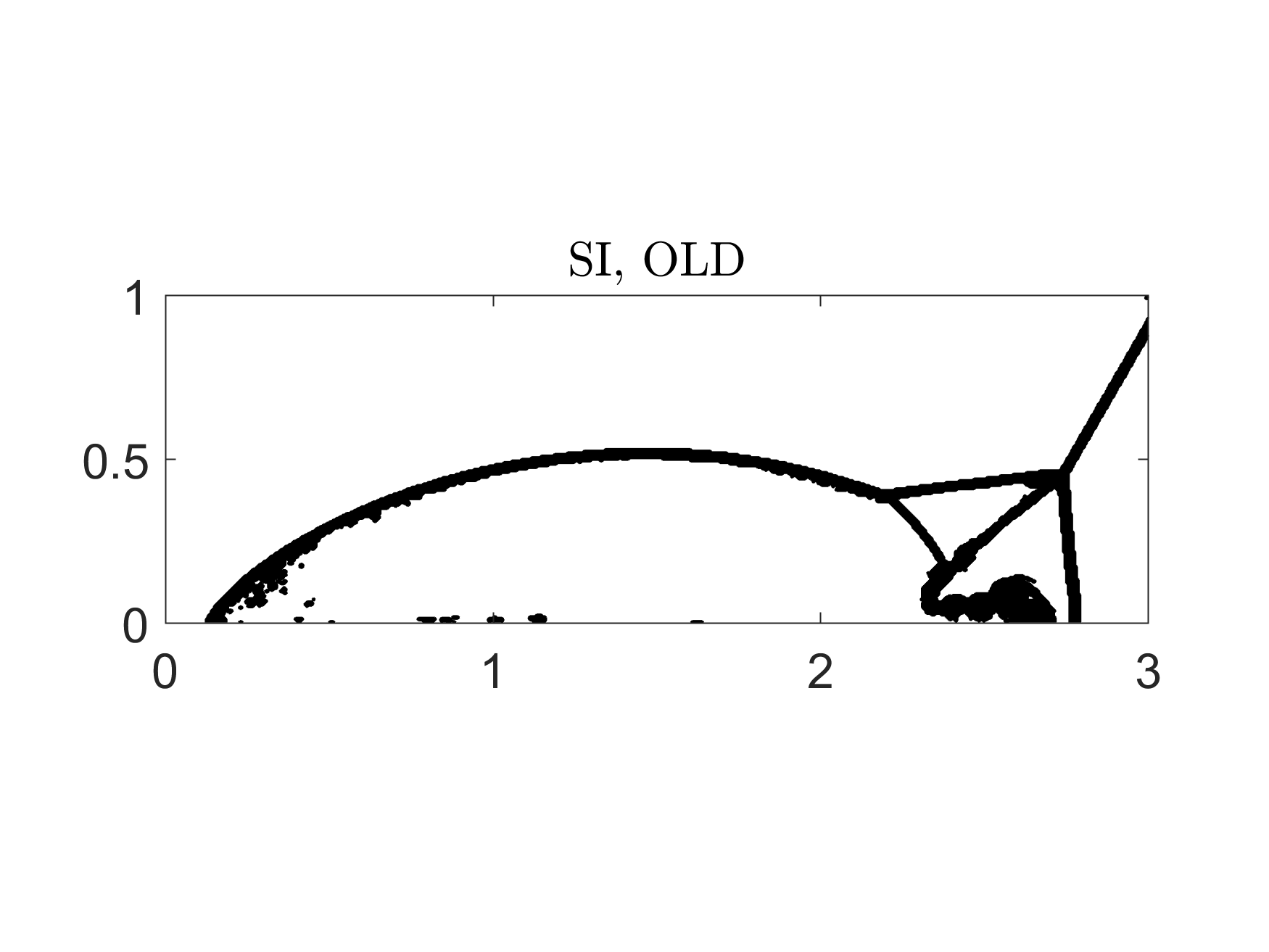}}
{\hspace{-3.4cm}\centerline{\includegraphics[trim=0.7cm 3.2cm 1.2cm 2.6cm, clip, width=7.3cm]{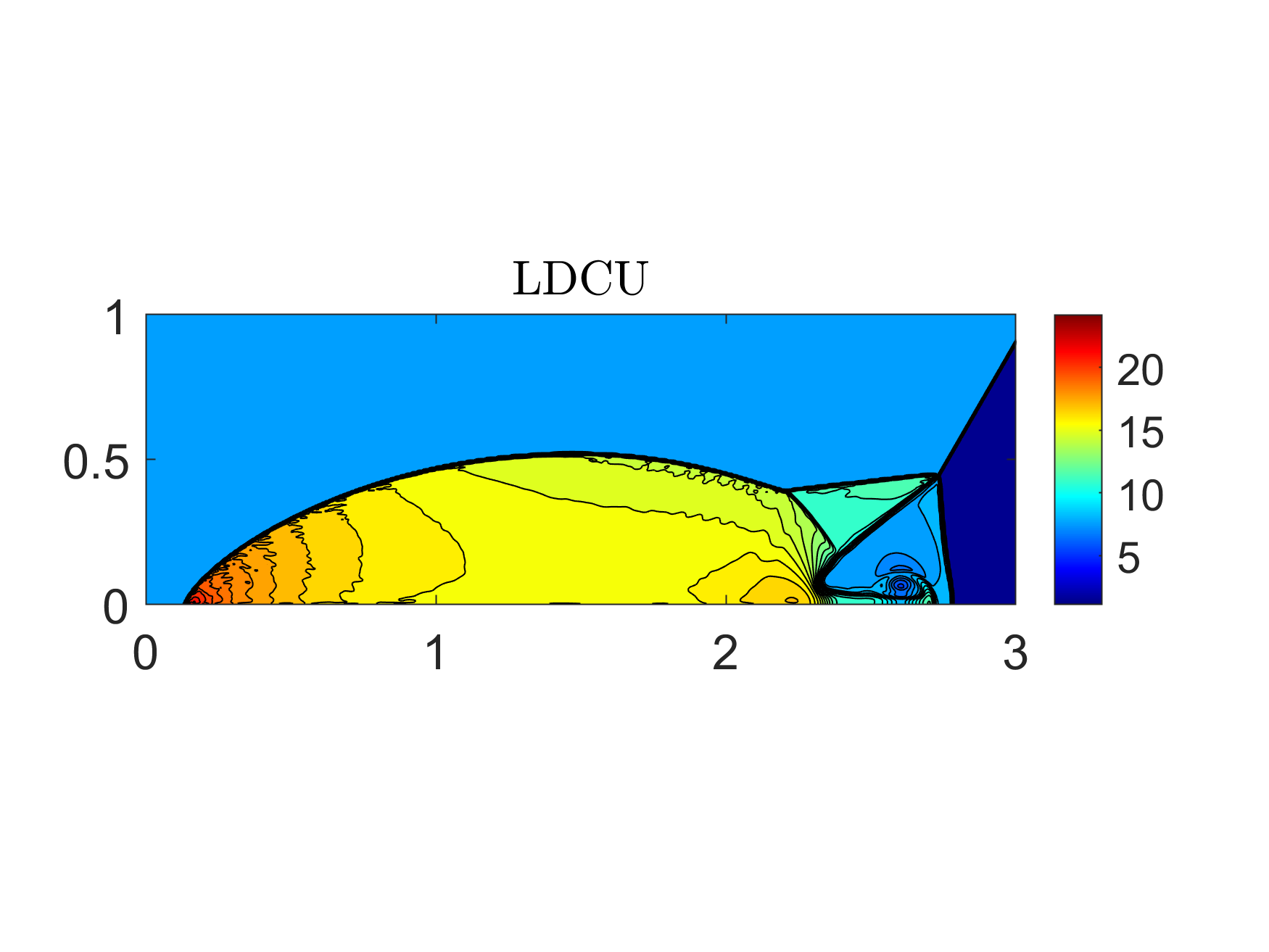}}}
\caption{\sf Example 8: Density $\rho$ computed by the NEW, OLD, and second-order LDCU schemes (left column) and the 
corresponding ``rough'' areas for the adaptive NEW and OLD schemes (right column).\label{fig8aa}}
\end{figure}

We compute the numerical solution until the final time $t=0.2$ by the studied NEW (with the adaption constants $\texttt{C}_1=0.004$ and
$\texttt{C}_2=0.4$) and OLD (with the adaption constant $\texttt{C}_1=0.01$) schemes as well as by the second-order LDCU scheme implemented
with the Minmod2 limiter. The results obtained on a uniform mesh with $\dx=\dy=1/200$ are in Figure \ref{fig8aa} (left column). As one can
see, both the NEW and OLD schemes clearly outperform the second-order LDCU scheme, and that the NEW scheme solution is slightly sharper and
substantially less oscillatory than the one obtained by the OLD scheme. At the same time, the ``rough'' areas detected by the NEW and OLD
schemes are quite similar; see Figure \ref{fig8aa} (right column).

\paragraph{Example 9---RT Instability.} In the last example taken from \cite{Shi03}, we investigate the RT instability, which is a physical
phenomenon occurring when a layer of heavier fluid is placed on top of a layer of lighter fluid. To this end, we first modify the 2-D Euler
equations of gas dynamics \eref{1.2}, \eref{3.3}--\eref{3.4} by adding the gravitational source terms acting in the positive direction of
the $y$-axis into the RHS of the system:
\begin{equation*}
\begin{aligned}
&\rho_t+(\rho u)_x+(\rho v)_y=0,\\
&(\rho u)_t+(\rho u^2 +p)_x+(\rho uv)_y=0,\\
&(\rho v)_t+(\rho uv)_x+(\rho v^2+p)_y=\rho,\\
&E_t+\left[u(E+p)\right]_x+\left[v(E+p)\right]_y=\rho v,
\end{aligned}
\end{equation*}
and then use the following initial conditions:
\begin{equation*}
(\rho,u,v,p)\Big|_{(x,y,0)}=
\begin{cases}(2,0,-0.025c\cos(8\pi x),2y+1),&y<0.5,\\(1,0,-0.025c\cos(8\pi x),y+1.5),&\mbox{otherwise},\end{cases}
\end{equation*}
where $c:=\sqrt{\gamma p/\rho}$ is the speed of sound. The solid wall boundary conditions are imposed at $x=0$ and $x=0.25$, and the
following Dirichlet boundary conditions are specified at the top and bottom boundaries:
$$
(\rho,u,v,p)(x,1,t)=(1,0,0,2.5),\quad(\rho,u,v,p)(x,0,t)=(2,0,0,1).
$$

\begin{figure}[ht!]
\centerline{\includegraphics[trim=5.cm 4.4cm 2.9cm 3.6cm, clip, width=14.cm]{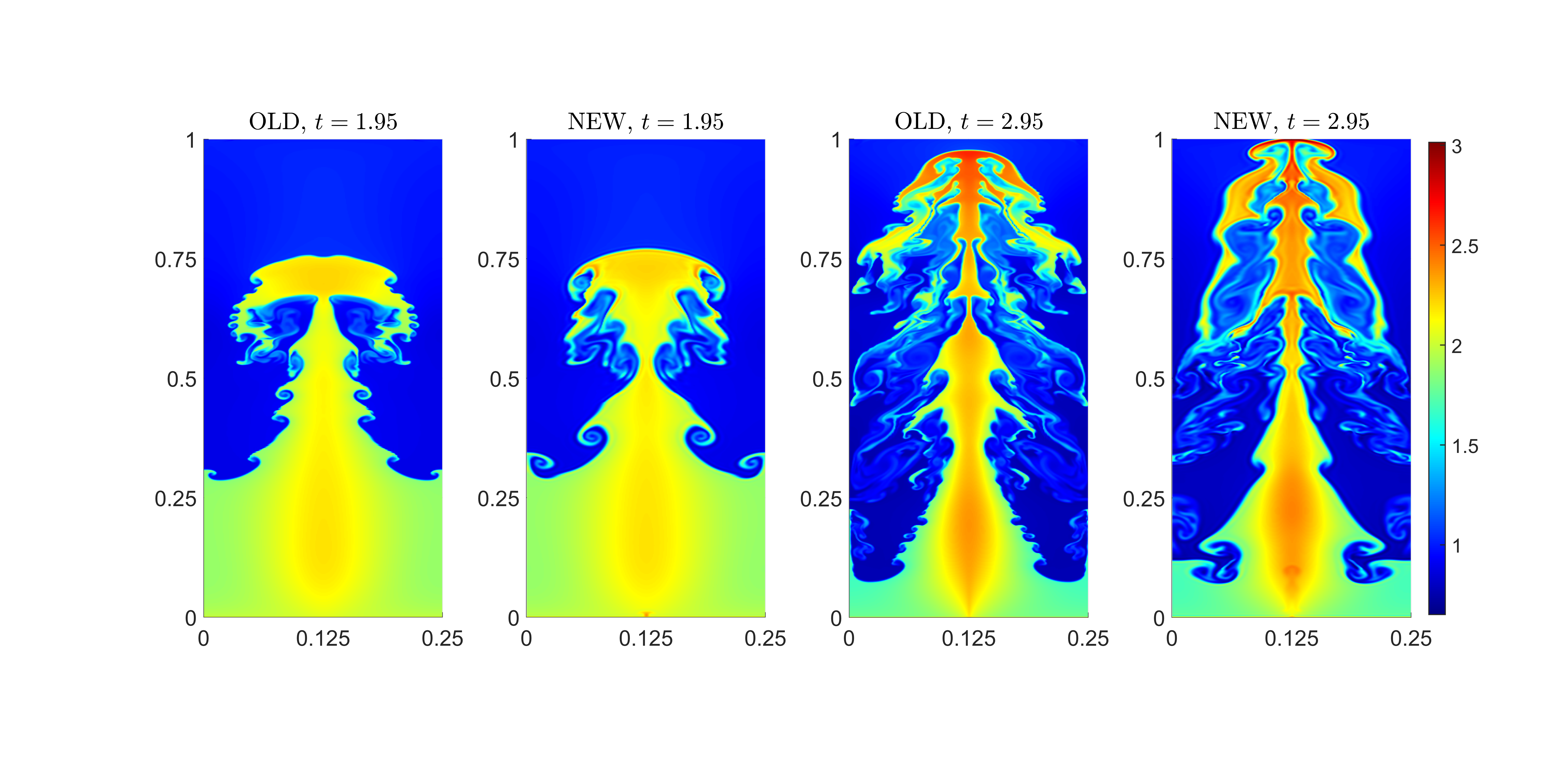}}
\caption{\sf Example 9: Density $\rho$ computed by the OLD and NEW adaptive schemes.\label{fig10a}}
\end{figure}

\begin{figure}[ht!] 
\centerline{\hspace{-0.3cm}\includegraphics[trim=4.3cm 0.5cm 4.4cm 0.2cm, clip, width=3.5cm]{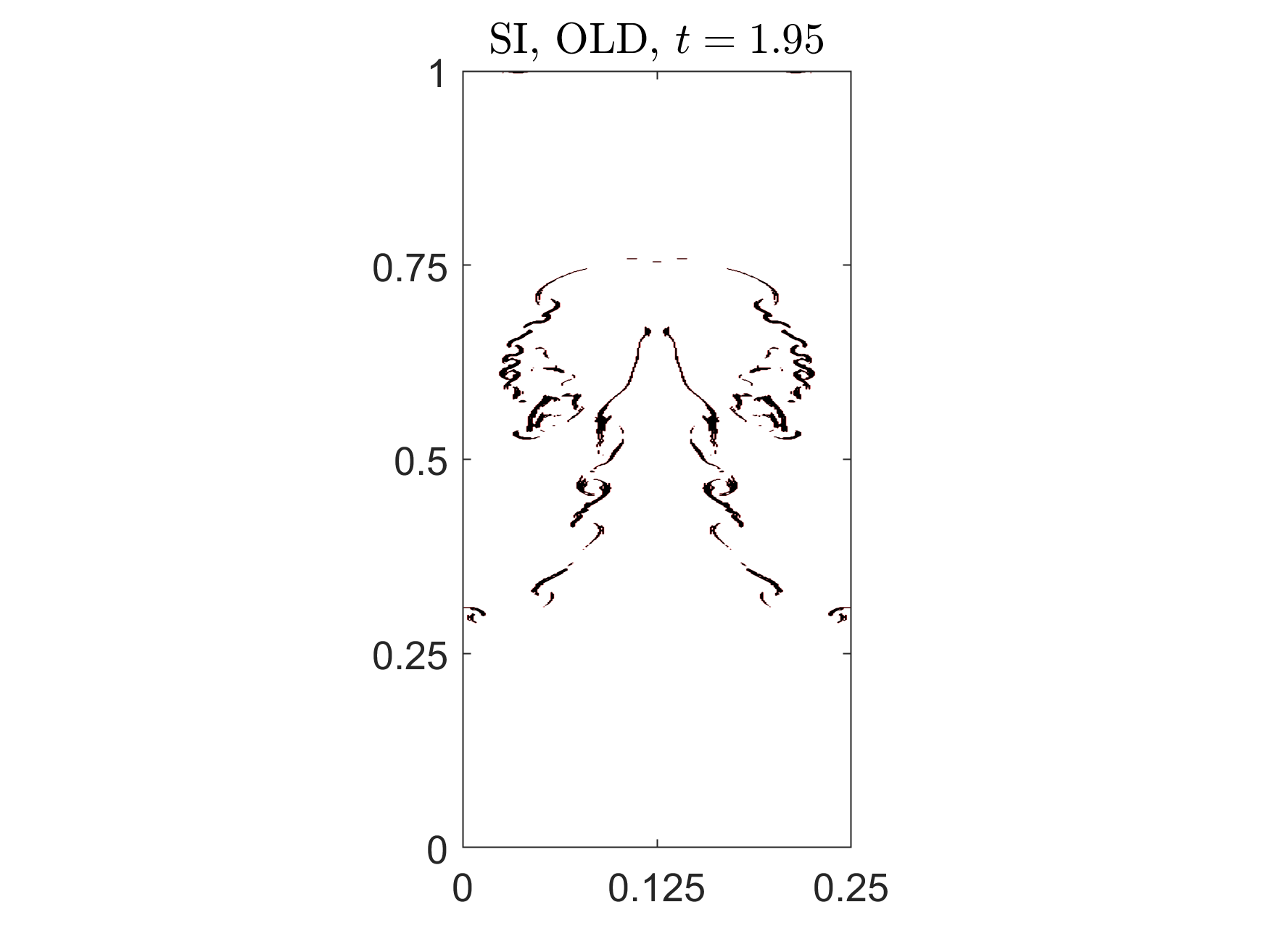}
                           \includegraphics[trim=4.3cm 0.5cm 4.4cm 0.2cm, clip, width=3.5cm]{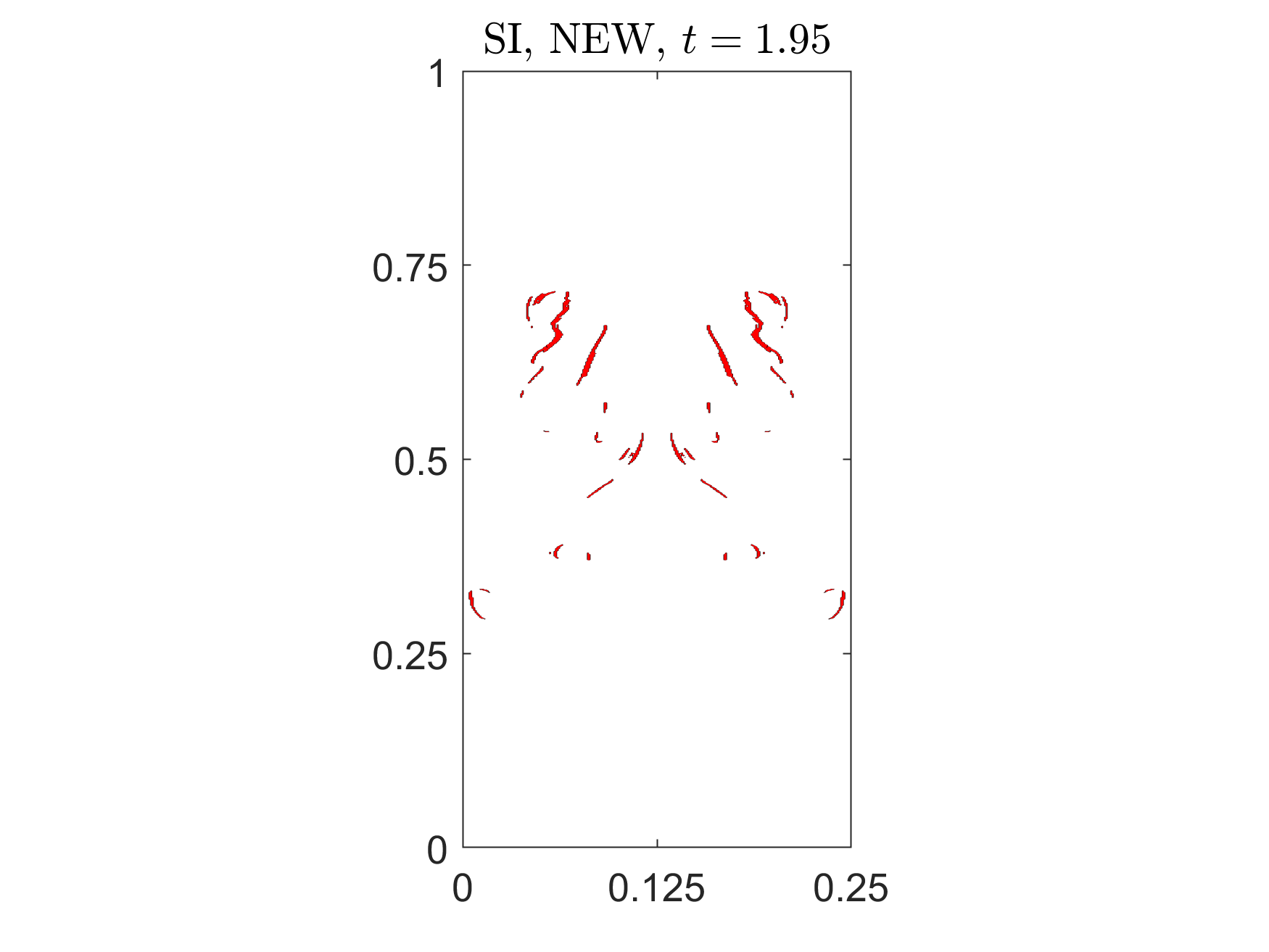}
                           \includegraphics[trim=4.3cm 0.5cm 4.4cm 0.2cm, clip, width=3.5cm]{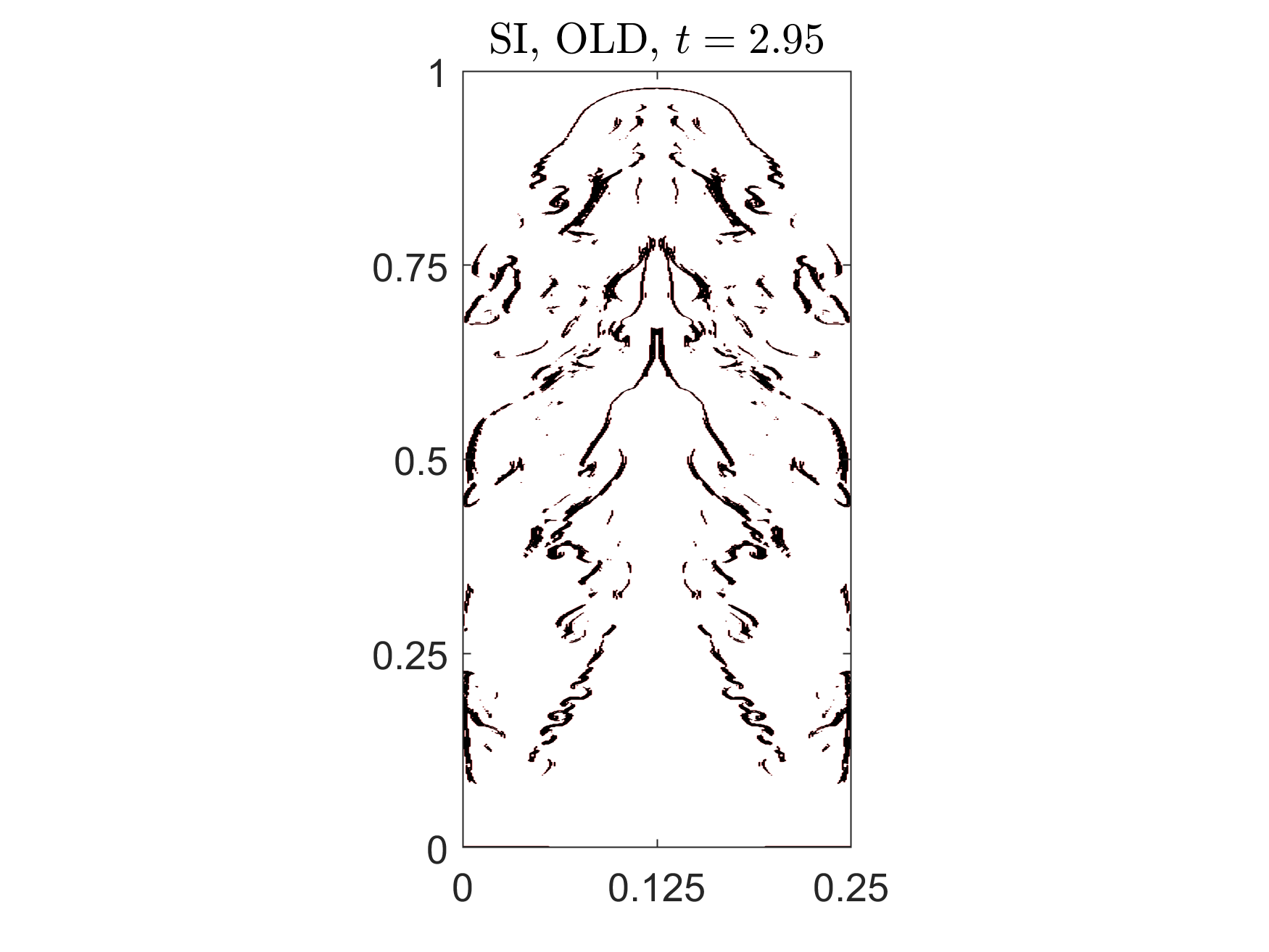}
                           \includegraphics[trim=4.3cm 0.5cm 4.4cm 0.2cm, clip, width=3.5cm]{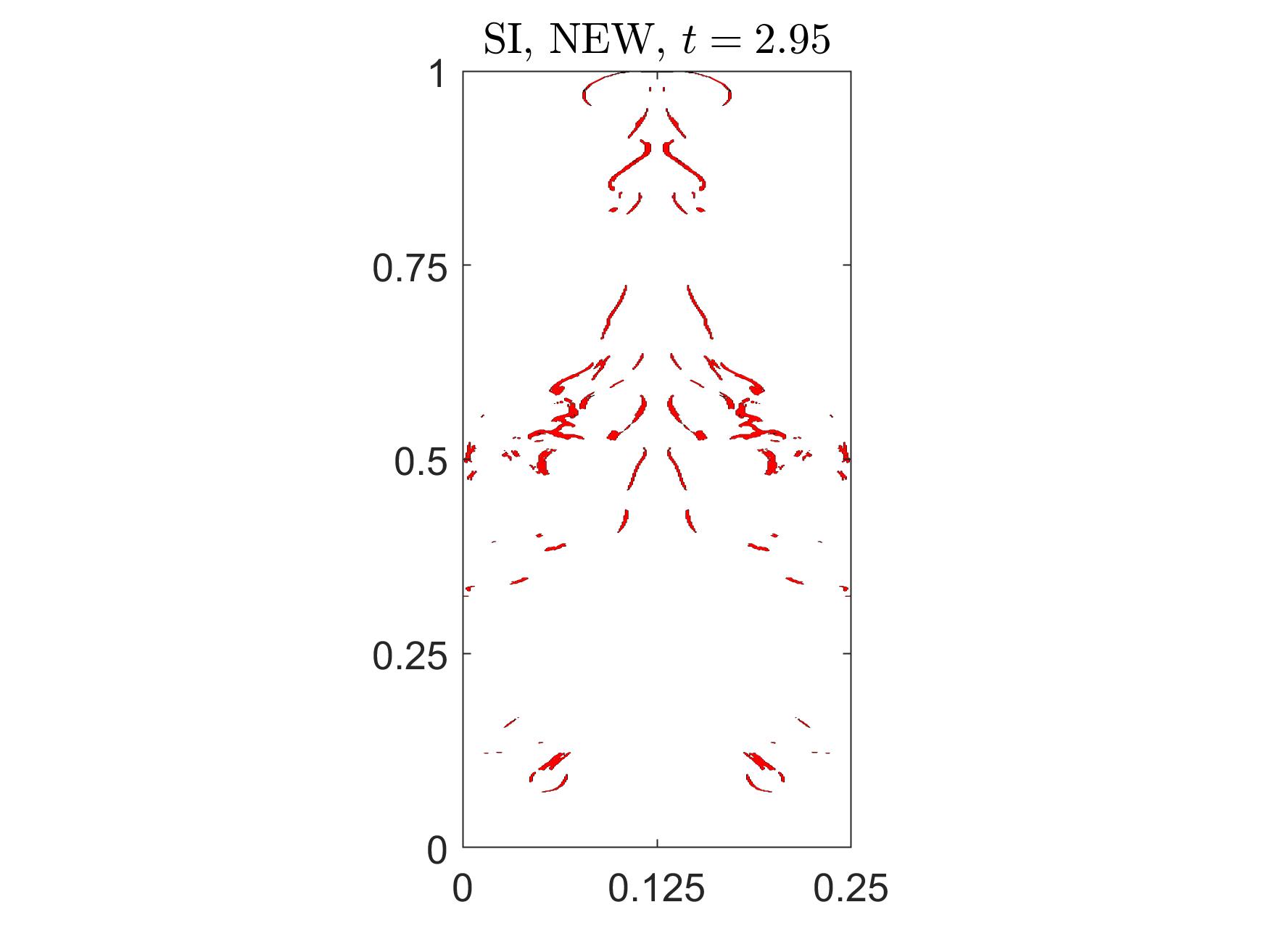}}
\caption{\sf Example 9: The ``rough'' areas detected by the OLD (left) and NEW (right) schemes at the final time-step.\label{fig10b}}
\end{figure}

We compute the numerical solution until the final time $t=2.95$ by the NEW (with the adaption constants $\texttt{C}_1=0.08$ and
$\texttt{C}_2=0.008$) and OLD (with the adaption constant $\texttt{C}_1=0.08$) schemes on a uniform mesh with $\dx=\dy=1/1024$ in the
computational domain $[0,0.25]\times[0,1]$. The obtained numerical results are presented in Figure \ref{fig10a}, where one can see that the
NEW scheme is less dissipative than the OLD scheme. As in the previous examples, we also show the ``rough'' areas detected by the studied
schemes; see Figure \ref{fig10b}, where one can see that in the NEW scheme, the overcompressive limiter is used in smaller areas compared
with the OLD scheme.

\section{Conclusion}\label{sec6}
We have introduced novel adaptive schemes for hyperbolic systems of conservation laws. The new schemes are based on a scheme adaption
strategy recently introduced in \cite{CKM2025}. Like any adaptive scheme, the proposed one utilizes a smoothness indicator (SI) to
automatically detect ``rough'' and smooth parts of the computed solution at each time-step. We then use the same SI to identify vicinities
of contact discontinuities within the detected ``rough'' areas. In the case of Euler equations of gas dynamics, this is achieved by applying
the SI to the density field, which is discontinuous at contact waves, and to the pressure, which is continuous there.

After splitting the computational domain into the corresponding three areas---``contact'', ``rough'', and ``smooth''---we apply different
schemes inside each of them. In the ``contact'' areas, we use the low-dissipation central-upwind (LDCU) scheme with the overcompressive SBM
limiter, which helps to extremely sharply resolve contact discontinuities. In the rest of the ``rough'' areas, we use the LDCU scheme with
the dissipative Minmod2 limiter, which helps to both capture shock waves in a non-oscillatory manner and prevent appearance of artificial
staircase-like structures, which may be developed when the overcompressive limiter is used. Finally, in the ``smooth'' area, we use a
recently proposed very simple quasi-linear fifth-order finite-difference scheme, which helps to extremely accurately capture smooth parts of
the computed solution.

The new adaptive schemes have been tested on a number of one- and two-dimensional numerical examples, which clearly demonstrate the ability
of the proposed schemes to achieve very high resolution, while being robust. We also show that the new schemes clearly outperform the
adaptive schemes from \cite{CKM2025}, while being substantially more efficient thanks to a very low computational cost of the quasi-linear
fifth-order schemes.

\noindent{\bf Declarations:\,\,}

\paragraph{Funding.} The work of S. Chu was funded by the DFG--SPP 2183: Eigenschaftsgeregelte Umformprozesse with the Project(s)
HE5386/19-2,19-3 Entwicklung eines flexiblen isothermen Reckschmiedeprozesses f\"ur die eigenschaftsgeregelte Herstellung von
Turbinenschaufeln aus Hochtemperaturwerkstoffen (424334423) and by the Deutsche Forschungsgemeinschaft (DFG, German Research
Foundation)--SPP 2410 Hyperbolic Balance Laws in Fluid Mechanics: Complexity, Scales, Randomness (CoScaRa) within the Project(s) HE5386/27-1
(Zuf\"allige kompressible Euler Gleichungen: Numerik und ihre Analysis, 525853336). The work of A. Kurganov was supported in part by NSFC
grants 12171226 and W2431004.

\paragraph{Conflicts of interest.} On behalf of all authors, the corresponding author states that there is no conflict of interest.

\paragraph{Data and software availability.} The data that support the findings of this study and FORTRAN codes developed by the authors and
used to obtain all of the presented numerical results are available from the corresponding author upon reasonable request.

\appendix
\section{1-D LCD-Based Piecewise Linear Reconstruction}\label{sec211}
In this appendix, we briefly describe a family of SBM limiters introduced in \cite{Lie03}, which are applied to the local characteristic
variables to compute the one-sided point values $\mU^\pm_\jph$. To this end, we first introduce the matrices
$\widehat A_\jph=A\big(\widehat{\mU}_\jph\big)$, where $\widehat{\mU}_\jph=\big(\,\xbar\mU_j+\,\xbar\mU_{j+1}\big)/2$, and compute the
matrices $R_\jph$ and $R^{-1}_\jph$ such that $R^{-1}_\jph\widehat A_\jph R_\jph$ are diagonal. We then introduce the local characteristic
variables $\bm\Gamma$ in the neighborhood of $x=x_\jph$:
$$
\bm\Gamma_\ell:=R^{-1}_\jph\,\xbar\mU_{j+\ell},\quad\ell=-1,0,1,2.
$$
Equipped with the values $\bm\Gamma_{-1}$, $\bm\Gamma_0$, $\bm\Gamma_1$, and $\bm\Gamma_2$, we compute the slopes
\begin{equation}
(\bm\Gamma_x)_0=\phi^{\rm SBM}_{\theta,\tau}\left(\frac{\bm\Gamma_1-\bm\Gamma_0}{\bm\Gamma_0-\bm\Gamma_{-1}}\right)
\frac{\bm\Gamma_0-\bm\Gamma_{-1}}{\dx}\quad\mbox{and}\quad(\bm\Gamma_x)_1=\phi^{\rm SBM}_{\theta,\tau}\left(\frac{\bm\Gamma_2-\bm\Gamma_1}
{\bm\Gamma_1-\bm\Gamma_0}\right)\frac{\bm\Gamma_1-\bm\Gamma_0}{\dx},
\label{2.4}
\end{equation}
where the two-parameter SBM function
\begin{equation}
\phi^{\rm SBM}_{\theta,\tau}(r):=\begin{cases}0,&r<0,\\\min\{r\theta,1+\tau(r-1)\},&0<r\le1,\\
r\phi^{\rm SBM}_{\theta,\tau}(\frac{1}{r}),&\text{otherwise,}\end{cases}
\label{2.6}
\end{equation}
is applied in the component-wise manner. Equipped with \eref{2.4}, we evaluate
\begin{equation}
\bm\Gamma^-_\hf=\bm\Gamma_0+\frac{\dx}{2}(\bm\Gamma_x)_0\quad\mbox{and}\quad\bm\Gamma^+_\hf=\bm\Gamma_1-\frac{\dx}{2}(\bm\Gamma_x)_1,
\label{2.7ff}
\end{equation}
and then obtain the corresponding point values of $\mU$ by $\mU^\pm_\jph=R_\jph\bm\Gamma^\pm_\hf$.

\noindent{\bf Remark A.1\,\,}
The parameters $\theta\in[1,2]$ and $\tau$ in \eref{2.6} can be used to adjust the amount of numerical dissipation present in the resulting
scheme. Generally, larger values of $\theta$ lead to sharper but potentially more oscillatory reconstructions. In all of the numerical
examples presented in this paper, we have taken $\theta=2$.

According to \cite{Lie03}, the behavior of the SBM limiter depends on $\tau$ in the following manner:

\smallskip
\noindent
$\bullet$ When $\tau\ge0.5$, the limiter is {\em dissipative}, which often leads to strong smearing of contact discontinuities over time;

\smallskip
\noindent
$\bullet$ For $0\le\tau<0.5$, the limiter is {\em compressive}, sharply resolving contact waves within few grid points; however, smooth
extrema may be slightly compressed, introducing small kinks in otherwise continuous profiles;

\smallskip
\noindent
$\bullet$ If $\tau<0$, the limiter is {\em overcompressive}, preserving contact discontinuities sharply over long times but overcompressing
smooth regions, possibly producing artificial ${\cal O}(1)$ jump discontinuities there.

\noindent{\bf Remark A.2\,\,}
The matrices $R_\jph$ and $R^{-1}_\jph$ in the case of the Euler equations of gas dynamics can be found in \cite[Appendix A]{CCHKL_22}.

\section{2-D LCD-Based Piecewise Linear Reconstruction}\label{appa}
In this appendix, we describe how to reconstruct the one-sided point values $\mU^\pm_{\jph,k}$ (the point values $\mU^\pm_{j,\kph}$ can be
computed in a similar manner and we omit the details for the sake of brevity). To this end, as in the 1-D case, we first introduce the local
characteristic variables in the neighborhood of $(x,y)=(x_\jph,y_k)$:
\begin{equation*}
\bm\Gamma_\ell:=R^{-1}_{\jph,k}\xbar\mU_{j+\ell,k},\quad\ell=-1,0,1,2,
\end{equation*}
where the matrix $R_{\jph,k}$ is such that $R^{-1}_{\jph,k}\widehat A_{\jph,k}R_{\jph,k}$ is diagonal and a locally linearized Jacobian is
$\widehat A_{\jph,k}:=A\big((\,\xbar\mU_{j,k}+\xbar\mU_{j+1,k})/2\big)$. We then compute $\bm\Gamma^\pm_\hf$ using \eref{2.4}--\eref{2.7ff},
and finally obtain the corresponding point values of $\mU$ by $\mU^\pm_{\jph,k}=R_{\jph,k}\bm\Gamma^\pm_\hf$.

\noindent{\bf Remark B.1\,\,}
The matrices $R_{\jph,k}$ and $R^{-1}_{\jph,k}$ for the 2-D Euler equations of gas dynamics \eref{1.2}, \eref{3.3}--\eref{3.4} can be found
in \cite[Appendix C]{CCHKL_22}.


\end{document}